\documentclass[a4paper,reqno,makeidx]{amsart}
\RequirePackage{ifthen}
\provideboolean{useutopia}
\setboolean{useutopia}{false}
\provideboolean{usesvjour}
\setboolean{usesvjour}{false}
\ifthenelse{\boolean{useutopia}}{%
  \RequirePackage{ifthen}
  \RequirePackage{iftex}
  \provideboolean{usemathrsfs}
  \provideboolean{useutopia}
  \setboolean{useutopia}{false}
  \ifXeTeX
    \RequirePackage[american,british]{babel}
    \RequirePackage{xltxtra}
    \DeclareSymbolFont{usualmathcal}{OMS}{cmsy}{m}{n}
    \DeclareSymbolFontAlphabet{\mathcalbf}{usualmathcal}
    \ifthenelse{\boolean{useutopia}}{
      \RequirePackage[utopia,euro]{mathdesign}
      \RequirePackage[OMLmathbf,OMLmathsfit]{isomath}
    }{
      \RequirePackage{amssymb}    
      \RequirePackage{bbold}    
    }
    \RequirePackage{mathrsfs} 
    \setboolean{usemathrsfs}{true}
  \else
    \RequirePackage[utf8]{inputenc}
    \ifthenelse{\boolean{useutopia}}{%
      \RequirePackage[utopia,euro]{mathdesign}
      \RequirePackage[OMLmathrm,OMLmathbf,OMLmathsfit]{isomath}
      \DeclareSymbolFont{usualmathcal}{OMS}{cmsy}{m}{n}
      \DeclareSymbolFontAlphabet{\mathcalbf}{usualmathcal}
      \RequirePackage{stmaryrd} 
      \providecommand{\diracdelta}[1][]{\ensuremath{\deltaup_{#1}}}
      
      \providecommand{\lap}{\ensuremath{\Deltaup}}
      \providecommand{\measure}[1]{\ensuremath{\mathcalbf{\uppercase{#1}}}}

    }{
      \RequirePackage{amssymb}    
      \RequirePackage{mathrsfs} 
      \RequirePackage{bbold}    
      \RequirePackage{stmaryrd} 
      \setboolean{usemathrsfs}{true}
      \providecommand{\mathcalbf}{\mathcal}
    }
    \RequirePackage[american,british]{babel}
  \fi
  \RequirePackage{mathtools}
  \RequirePackage{nicefrac}
  \RequirePackage{stmaryrd} 
  \RequirePackage{amsthm}
  \RequirePackage{enumerate}
  \RequirePackage{xspace}
  \RequirePackage{verbatim}
  \RequirePackage{fancyvrb}
  \RequirePackage{hyperref}
  \RequirePackage{listings}
  \RequirePackage{xcolor}
  \RequirePackage{epsfig}
  \RequirePackage{graphicx}
  \RequirePackage{booktabs}
  \RequirePackage{tikz}
  \usetikzlibrary{calc}
  \provideboolean{usebeamer}
  \provideboolean{isthesis}
  \provideboolean{isamsltex}
  \provideboolean{issiamltex}
  
}{%
  \RequirePackage{ifthen}
  \RequirePackage{iftex}
  \provideboolean{usemathrsfs}
  \provideboolean{useutopia}
  \setboolean{useutopia}{false}
  \ifXeTeX
    \RequirePackage[american,british]{babel}
    \RequirePackage{xltxtra}
    \DeclareSymbolFont{usualmathcal}{OMS}{cmsy}{m}{n}
    \DeclareSymbolFontAlphabet{\mathcalbf}{usualmathcal}
    \ifthenelse{\boolean{useutopia}}{
      \RequirePackage[utopia,euro]{mathdesign}
      \RequirePackage[OMLmathbf,OMLmathsfit]{isomath}
    }{
      \RequirePackage{amssymb}    
      \RequirePackage{bbold}    
    }
    \RequirePackage{mathrsfs} 
    \setboolean{usemathrsfs}{true}
  \else
    \RequirePackage[utf8]{inputenc}
    \ifthenelse{\boolean{useutopia}}{%
      \RequirePackage[utopia,euro]{mathdesign}
      \RequirePackage[OMLmathrm,OMLmathbf,OMLmathsfit]{isomath}
      \DeclareSymbolFont{usualmathcal}{OMS}{cmsy}{m}{n}
      \DeclareSymbolFontAlphabet{\mathcalbf}{usualmathcal}
      \RequirePackage{stmaryrd} 
      \providecommand{\diracdelta}[1][]{\ensuremath{\deltaup_{#1}}}
      
      \providecommand{\lap}{\ensuremath{\Deltaup}}
      \providecommand{\measure}[1]{\ensuremath{\mathcalbf{\uppercase{#1}}}}

    }{
      \RequirePackage{amssymb}    
      \RequirePackage{mathrsfs} 
      \RequirePackage{bbold}    
      \RequirePackage{stmaryrd} 
      \setboolean{usemathrsfs}{true}
      \providecommand{\mathcalbf}{\mathcal}
    }
    \RequirePackage[american,british]{babel}
  \fi
  \RequirePackage{mathtools}
  \RequirePackage{nicefrac}
  \RequirePackage{stmaryrd} 
  \RequirePackage{amsthm}
  \RequirePackage{enumerate}
  \RequirePackage{xspace}
  \RequirePackage{verbatim}
  \RequirePackage{fancyvrb}
  \RequirePackage{hyperref}
  \RequirePackage{listings}
  \RequirePackage{xcolor}
  \RequirePackage{epsfig}
  \RequirePackage{graphicx}
  \RequirePackage{booktabs}
  \RequirePackage{tikz}
  \usetikzlibrary{calc}
  \provideboolean{usebeamer}
  \provideboolean{isthesis}
  \provideboolean{isamsltex}
  \provideboolean{issiamltex}
  
}
\colorlet{a}{magenta}
\colorlet{b}{green!75!blue}
\colorlet{c}{yellow!87.5!red}
\colorlet{d}{cyan}
\colorlet{e}{red}
\colorlet{f}{blue}
\colorlet{g}{white}
\colorlet{h}{black!50}
\colorlet{i}{black}
\colorlet{j}{black!75}
 
\providecommand{\linkedurl}[1]{\url{#1}}
\providecommand{\linkedemail}[1]{\href{mailto:#1}{#1}}

\providecommand{\email}[1]{{\linkedemail{#1}}}

\providecommand{\Ignore}[1]{}
\providecommand{\ignore}[1]{}
\providecommand{\freeze}[1]{}
\providecommand{\crossout}[1]{{\textcolor{i!20}{#1}}}
\providecommand{\highlightcolor}{a}
\providecommand{\highlight}[1]{{\color{\highlightcolor}#1}}

\providecommand{\memo}[1]{%
  \ensuremath{%
    \framebox{\tiny\textbf{\kern-2pt\textsf{#1}}\kern-2pt}%
  }%
  \xspace}

\RequirePackage{alphalph}
\provideboolean{shownotes}
\setboolean{shownotes}{true}
\newcounter{margnote}[page]
\providecommand{\mgcolor}{a}
\providecommand{\mgcolorset}[1]{\renewcommand{\mgcolor}{\alphalph{#1}}}
\providecommand{\mgcolorsetbycounter}[1]{%
  \ifthenelse{\value{#1}<11}{%
    \renewcommand{\mgcolor}{\alph{#1}}%
  }{%
    \renewcommand{\mgcolor}{a}}%
}
\providecommand{\mgcolormake}{\mgcolorsetbycounter{margnote}}
\providecommand{\mgcolorstepby}[1]{
  \setcounter{tmpcounter}{\value{margnote}}%
  \addtocounter{tmpcounter}{#1}%
  \mgcolorsetbycounter{tmpcounter}%
}%
\providecommand{\margnotecolor}{%
  \ifthenelse{\value{margnote}=0}{%
    \mgcolorset{10}
  }{%
    \ifthenelse{\value{margnote}<7}{%
      \mgcolormake%
    }{%
      \ifthenelse{\value{margnote}=7}{\mgcolorset{10}}{%
        \ifthenelse{\value{margnote}<11}{\mgcolormake}{%
          \ifthenelse{\value{margnote}<17}{\mgcolorstepby{-10}}{%
            \mgcolorset{10}%
          }%
        }%
      }%
    }%
  }%
}%
\providecommand{\margnotemark}{{\colorbox{\mgcolor}{\tiny\color{g}\upshape\texttt{\arabic{page}.\arabic{margnote}}}}}
\providecommand{\margnote}[2][]{%
  \ifthenelse{%
    \boolean{shownotes}%
  }{%
    \stepcounter{margnote}%
    \margnotecolor%
    \margnotemark %
    \marginpar{%
      \color{\mgcolor}%
      \texttt{%
        \begin{minipage}{2cm}%
          \raggedright\tiny%
          \margnotemark%
          #2%
          \\
          {\ifx|#1|{}\else{ - #1}\fi}%
        \end{minipage}%
      }%
    }%
  }{%
  }%
}%
\providecommand{\mathnote}[2][]{%
  \ifthenelse{%
    \boolean{shownotes}%
  }{%
    \stepcounter{margnote}%
    \margnotecolor%
    \text{%
      \colorbox{\mgcolor}{%
        \color{g}%
        \texttt{%
          \tiny%
              \margnotemark: 
              \ifx|#1|{}\else{#1:}\fi%
              #2%
        }%
      }%
    }%
  }{%
  }%
}%
\providecommand{\textnote}[2][]{%
  \ifthenelse{%
    \boolean{shownotes}%
  }{%
    \stepcounter{margnote}%
    \margnotecolor%
    \ \\
    \text{%
      \colorbox{\mgcolor}{%
        \begin{minipage}{.9\textwidth}
        \color{g}%
        \texttt{%
          \margnotemark: 
          \ifx|#1|{}\else{#1: }\fi%
          #2%
        }%
        \end{minipage}
      }%
    }%
  }{%
  }%
}%

\providecommand{\Todo}[1]{
  \ifthenelse{\boolean{shownotes}}{
    \begin{center}
    \begin{tikzpicture}
     \node[fill=a!17]{
       \begin{minipage}{\textwidth}
         \texttt{To do:}
         \\
         \texttt{\bfseries{\small #1}}
       \end{minipage}
     };
    \end{tikzpicture}
    \end{center}
  }{}}
\provideboolean{showrevisions}
\setboolean{showrevisions}{true}
\newcommand{\revisionsheader}{***\newline\Warning{the following part is under revision}}
\newcommand{\revisionsfooter}{***\newline\Warning{end of part under revision}}

\providecommand{\Warning}[1]{    
  \begin{tikzpicture}
    \node[fill=a!27]{
      \begin{minipage}{\textwidth}
        \texttt{\bfseries{\small Warning: #1}}
      \end{minipage}
    };
  \end{tikzpicture}
}

\provideboolean{showcomments}
\providecommand{\margincomment}[1]{
\ifthenelse{\boolean{showcomments}}{\marginpar{\tiny #1}}{}
}
\provideboolean{showchanges}
\setboolean{showchanges}{false}
\providecommand{\changes}[2][]{
  \ifthenelse{\boolean{showchanges}}{{\ifx|#1|{}\else\margnote{#1}\fi\highlight{#2}}}{#2}}
\providecommand{\mathchanges}[2][]{
  \ifthenelse{\boolean{showchanges}}{{\ifx|#1|{}\else\mathnote{#1}\fi\highlight{#2}}}{#2}}
\newenvironment{LongChanges}{
  \ifthenelse{\boolean{showchanges}}{\color{\highlightcolor}}{}
}{}
\providecommand{\changefromto}[3][replace with]{%
  \ifthenelse{\boolean{showchanges}}{%
    {\crossout{#2}\margnote{#1}}{\highlight{#3}}}{%
    #3\xspace}%
}
\providecommand{\ChangePar}[3][]{
  \ifthenelse{\boolean{showchanges}}{
    {\par\textcolor{i!20}{#2}\ifx|#1|\else\margnote{#1}\fi}{\par\textcolor{a}{#3}}
  }{%
    \par #3%
  }%
}
\providecommand{\InsertPar}[1]{
  \ifthenelse{\boolean{showchanges}}
  {{\par$\mapsto$ \textcolor{blue}{#1}}}
  {\par #1}
}

\providecommand{\mathchangefromto}[3][]{\crossout{#2}\ifx|#1|\else\mathnote{#1}\fi\highlight{#3}}

\providecommand{\mathscript}
	   {\mathscr}

 \providecommand{\bbbold}{\mathbb}

 \providecommand{\rR}{\ensuremath{\bbbold R}\xspace}
 
 \providecommand{\rT}{\ensuremath{\bbbold T}\xspace}

\providecommand{\Ae}[1][]{\ensuremath{\ifx|#1|{\ }\else{\:#1\text{-}}\fi\text{almost everywhere }}\xspace}
\providecommand{\Aa}[1][]{\ensuremath{\text{ for }\ifx|#1|{}\else{\:#1\text{-}}\fi\text{almost all }}}
\providecommand{\as}[1][]{\ensuremath{\ifx|#1|{\ }\else{#1\text{-}}\fi\text{almost surely}}\xspace}

 \providecommand{\reals}{\rR}

 \providecommand{\R}[1]{\reals^{#1}}
 
 \providecommand{\fieldmats}[3][F]{\csname#1\endcsname{#2\times#3}}

 \providecommand{\RO}[1][]{{\reals_{0+}\ifx|#1|{}\else^{#1}\fi}}
 \providecommand{\RP}[1][]{{\reals_+\ifx|#1|{}\else^{#1}\fi}}

 \providecommand{\torus}[1]{\rT\ifthenelse{\equal{#1}1}{}{^#1}}

 \providecommand{\one}{\ensuremath{\bbbold 1}\xspace}
 \providecommand{\charfun}[1]{\one_{#1}}
 \providecommand{\iverson}[1]{\one_{\qb{#1}}}

 \providecommand{\diracdelta}[1][]{\ensuremath{{\mathrm{\delta}}\ifx|#1|{}\else_{#1}\fi}}

 \providecommand{\pic}{\ensuremath{\mathrm\pi}}
 \providecommand{\pifracl}[2][]{\fracl{\ifx|#1|\else#1\fi\pic}{#2}}
 \providecommand{\pifrac}[2][]{\frac{\ifx|#1|\else#1\fi\pic}{#2}}

 \providecommand{\inner}{\cdot}
 \providecommand{\outerp}{\wedge}

 \providecommand{\W}{\ensuremath{\varOmega}\xspace}

 \providecommand{\ep}{\ensuremath{\varepsilon}\xspace}

 \providecommand{\qgroup}[1]{{#1}}
 
 \providecommand{\qp}[2][]{\ifx|#1|\left(\else\csname#1\endcsname(\fi{#2}\ifx|#1|\right)\else\csname#1\endcsname)\fi}

 \providecommand{\qpreg}[1]{\ensuremath{(#1)}}
 \providecommand{\qpbig}[1]{\qp[big]{#1}}
 \providecommand{\qpBig}[1]{\ensuremath{\Big(#1\Big)}}
 \providecommand{\qpbigg}[1]{\ensuremath{\bigg(\!#1\!\bigg)}}
 \providecommand{\qpBigg}[1]{\ensuremath{\Bigg(\!#1\!\Bigg)}}
 \providecommand{\qb}[2][]{\ifx|#1|\left[\else\csname#1\endcsname[\fi{#2}\ifx|#1|\right]\else\csname#1\endcsname]\fi}
 \providecommand{\qc}[2][]{\ifx|#1|\left\{\else\csname#1\endcsname\{\fi{#2}\ifx|#1|\right\}\else\csname#1\endcsname\}\fi}

 \providecommand{\qa}[1]{\ensuremath{\left\langle{#1}\right\rangle}}
 \providecommand{\qareg}[1]{\ensuremath{\langle#1\rangle}}
 \providecommand{\qabig}[1]{\ensuremath{\big\langle#1\big\rangle}}
 \providecommand{\qaBig}[1]{\ensuremath{\Big\langle#1\Big\rangle}}
 \providecommand{\qabigg}[1]{\ensuremath{\bigg\langle#1\bigg\rangle}}
 \providecommand{\qaBigg}[1]{\ensuremath{\Bigg\langle#1\Bigg\rangle}}

 \providecommand{\opinter}[2]{\ensuremath{\left(#1,#2\right)}\xspace}
 \providecommand{\clinter}[2]{\ensuremath{\left[#1,#2\right]}\xspace}

 \providecommand{\compowqp}[2]{\ensuremath{\qp{\!#2\!\!}^{\kern -.4em #1}\!}}
 
 \providecommand{\powqpreg}[2]{\ensuremath{%
     \qpreg{#2}^{\kern 0em\lower .1ex\hbox{\scriptsize $#1$}}\kern-.3em}}
 \providecommand{\powqpbig}[2]{\ensuremath{%
     \qpbig{#2}^{\kern -.2em\lower .3ex\hbox{\scriptsize $#1$}}\kern-.3em}}
 \providecommand{\powqpBig}[2]{\ensuremath{%
     \qpBig{#2}^{\kern -.2em\lower .3ex\hbox{\scriptsize $#1$}}\kern-.3em}}
 \providecommand{\powqpbigg}[2]{\ensuremath{%
     \qpbigg{#2}^{\kern -.2em\lower .3ex\hbox{\scriptsize $#1$}}\kern-.3em}}
 \providecommand{\powqpBigg}[2]{\ensuremath{%
     \qpBigg{#2}^{\kern -.2em\lower .3ex\hbox{\scriptsize $#1$}}}}
 \providecommand{\powp}[3][]{#3\ifx|#1|^{#2}\else{#1}^{#2}\fi}
 \providecommand{\pow}[2][]{\ifx|#1|\operatorname{pow}^{#2}\else\powp{#2}{#1}\fi}

 \providecommand{\powabs}[2]{\powp{#1}{\abs{#2}}}

 \providecommand{\norm}[2][]{\ifx|#1|\left|\else\csname#1\endcsname|\fi#2\ifx|#1|\right|\else\csname#1\endcsname|\fi}
 \providecommand{\normon}[2]{\norm{#1}_{#2}}

 \providecommand{\abs}[2][]{\ensuremath{\ifx|#1|{\left|#2\right|}\else{\csname#1\endcsname|{#2}\csname#1\endcsname|}\fi}}

 \providecommand{\Norm}[2][]{\ifx|#1|\left\|\else\csname#1\endcsname\|\fi{#2}\ifx|#1|\right\|\else\csname#1\endcsname\|\fi}

 \providecommand{\Normon}[2]{\Norm{#1}_{#2}}

 \providecommand{\normonsob}[4][]{\normon{#2}{\sob{#3}{#4}\if|#1|{}\else(#1)\fi}}
 \providecommand{\Normonsob}[4][]{\Normon{#2}{\sob{#3}{#4}\if|#1|{}\else(#1)\fi}}
 \providecommand{\Normonleb}[4][]{\Normon{#2}{\leb{#3}\if|#1|{}\else(#1)\fi}}
 \providecommand{\ltwop}[3][]{\ensuremath{\qa{#2,#3}\ifx|#1|\else_{#1}\fi}}
 \providecommand{\ltwopreg}[2]{\ensuremath{\qareg{#1,#2}\ifx|#1|\else_{#1}\fi}}
 \providecommand{\ltwopbig}[2]{\ensuremath{\qabig{#1,#2}\ifx|#1|\else_{#1}\fi}}
 \providecommand{\ltwopBig}[2]{\ensuremath{\qaBig{#1,#2}\ifx|#1|\else_{#1}\fi}}
 \providecommand{\ltwopbigg}[2]{\ensuremath{\qabigg{#1,#2}\ifx|#1|\else_{#1}\fi}}
 \providecommand{\ltwopBigg}[2]{\ensuremath{\qaBigg{#1,#2}\ifx|#1|\else_{#1}\fi}}
 
 \providecommand{\average}[2][]{{\qa{#2}\ifx|#1|\else_{#1}\fi}}
 
 \providecommand{\ensemble}[2]{\ensuremath{\left\{ #1:\;#2 \right\}}}

 \providecommand{\setof}[1]{{\qc{#1}}}

 \providecommand{\conditionalto}[1]{{\left|{#1}\right.}}
 
\providecommand{\measure}[1]{\ensuremath{\mathcalbf{\MakeUppercase{#1}}}}

\providecommand{\probmeasure}[2][]{{\measure{#2}}\ifx|#1|\else_{#1}\fi}
\providecommand{\Prob}{}
\renewcommand{\Prob}[1][]{\probmeasure[{#1}]{p}}

\providecommand{\randvars}[1][\Prob]{\operatorname{RV}\ifx|#1|{}\else{(#1)}\fi}
\providecommand{\discrandvars}[1][\Prob]{\operatorname{DRV}\ifx|#1|{}\else{({#1)}\fi}} 
\providecommand{\contrandvars}[1][\Prob]{\ensuremath{\operatorname{CDRV}\ifx|#1|{}\else(#1)\fi}} 
 \def\env@matrix{\hskip -\arraycolsep
  \let\@ifnextchar\new@ifnextchar
  \array{*\c@MaxMatrixCols c}}
 \makeatletter
 \renewcommand*\env@matrix[1][c]{\hskip -\arraycolsep
   \let\@ifnextchar\new@ifnextchar
   \array{*\c@MaxMatrixCols #1}}
 \makeatother
 \providecommand{\irow}[2]{#1_{#2}}
 \providecommand{\icol}[2]{#1^{#2}}

 \providecommand{\ijrowcol}[3]{\icol{\irow{#1}{#2}}{#3}}
 \providecommand{\entry}[1]{\qb{#1}}
 \providecommand{\vecentry}[2]{\irow{#1}{#2}}

 \providecommand{\vecof}[1]{\qp{#1}}
 
 \providecommand{\rowof}[1]{\qb{#1}}

 \providecommand{\getentryi}[2]{\irow{\entry{#1}}{#2}}
 
 \providecommand{\getvecentry}[2]{\getentryi{\vec #1}{#2}}

 \providecommand{\getrow}[2]{\irow{\entry{#1}}{#2}}

 \providecommand{\dismatof}[2][r]{\begin{bmatrix}[#1]#2\end{bmatrix}}

 \providecommand{\matentry}[3]{\ijrowcol{#1}{#2}{#3}}

 \providecommand{\block}[5]{\ijrowcol{#1}{\ifx#2#3{\rowof{#2}}\else\rowof{{#2}\dotsc{#3}}\fi}{\ifx#4#5{\rowof{#4}}\else\rowof{{#4}\dotsc{#5}}\fi}}
 \providecommand{\colblock}[3]{\getvecentry{#1}{\ifx#2#3{#2}\else\fromto{#2}{#3}\fi}}

 \providecommand{\dismatskeldots}[4]{
   \dismatof[c]{
     #1&\dotsc&#3
     \\
     \vdots & \ddots &\vdots
     \\
     #2&\dotsc&#4
   }
 }
 \providecommand{\dismatcommfromtofromto}[5]{
   \dismatskeldots{#1#2#4}{#1#3#4}{#1#2#5}{#1#3#5}
 }
 \providecommand{\dismatcustfromtofromto}[6][matentry]{
   \dismatcommfromtofromto{\csname#1\endcsname{#2}}#3#4#5#6
 }
 \providecommand{\dismatcustfromtofromto}[6][matentry]{
   \dismatskeldots{%
     \csname#1\endcsname{#2}{#3}{#4}%
   }{%
     \csname#1\endcsname{#2}{#3}{#6}%
   }{%
     \csname#1\endcsname{#2}{#5}{#4}%
   }{%
     \csname#1\endcsname{#2}{#5}{#6}%
   }%
 }%
 \providecommand{\dismatcustfromtofromto}[6][matentry]{
   \dismatof{
     \csname#1\endcsname{#2}{#3}{#4}&\dotsc&\csname#1\endcsname{#2}{#3}{#6}
     \\
     \vdots & \ddots &\vdots
     \\
     \csname#1\endcsname{#2}{#5}{#4}&\dotsc&\csname#1\endcsname{#2}{#5}{#6}
   }
 }

 \providecommand{\dissysaxbdotsnm}[5]{\begin{matrix}[r]%
     \matentry{#1}11\vecentry{#2}1&+\dotsb&+\matentry{#1}1{#5}\vecentry{#2}{#5}
     &
     =
     \ifx|#3|0\else{\vecentry {#3}1}\fi
     \\
     \dotsb
     \\
     \matentry{#1}{#4}1\vecentry{#2}1&+\dotsb&+\matentry{#1}{#4}{#5}\vecentry{#2}{#5}
     &
     =
     \ifx|#3|0\else{\vecentry {#3}{#4}}\fi
 \end{matrix}}
 
 \providecommand{\seqof}[1]{\qp{#1}}
 \providecommand{\seqs}[2]{\seqof{#1}_{#2}}
 \providecommand{\sets}[2]{\setof{#1}_{#2}}
 \providecommand{\seqi}[3][]{\seqs{#2_{#3}}{\ifx|#1|{#3}\else{{#3}\in{#1}}\fi}}

 \providecommand{\seti}[3][]{\sets{#2_{#3}}{\ifx|#1|_{#3}\else_{{#3}\in{#1}}\fi}}
 \providecommand{\sequ}[3][]{\seqs{#2^{#3}}{\ifx|#1|{#3}\else{{#3}\in{#1}}\fi}}
 \providecommand{\setu}[3][]{\sets{#2^{#3}}{\ifx|#1|{#3}\else{{#3}\in{#1}}\fi}}

 \providecommand{\limofat}[3][]{\ensuremath{\lim_{\ifx|#1|{}\else{#1\ni}\fi#3}{#2}}}
 \providecommand{\limsupofat}[3][]{\ensuremath{\limsup_{\ifx|#1|{}\else{#1\ni}\fi#3}{#2}}}
 \providecommand{\liminfofat}[3][]{\ensuremath{\liminf_{\ifx|#1|{}\else{#1\ni}\fi#3}{#2}}}

 \providecommand{\jump}[2][]{\ensuremath{\left\llbracket #2\right\rrbracket\ifx|#1|{}\else_{#1}\fi}}

 \providecommand{\fromto}[2]{\ensuremath{\setof{#1\dotsc#2}}}
 
 \providecommand{\d}{}
 \renewcommand{\d}[1][]{\ensuremath{\operatorname{d}\!\ifx|#1|\else{_{#1}}\fi}}
 
 \providecommand{\ds}[1][]{\d{\measure S}}
 \providecommand{\D}[1][]{\ensuremath{\operatorname{D}\!\ifx|#1|\else{_{#1}}\fi}}
 
\providecommand{\registered}%
{\ensuremath{^\text{\textregistered}}}

\providecommand{\tand}{\ensuremath{\text{ and }}}

\providecommand{\const}{\operatorname{constant}}     
\providecommand{\constant}[1]{\ensuremath{C_{#1}}}
\providecommand{\constext}[2][]{\constant{\textup{#2}{\ifx|#1|{}\else{,\ensuremath{#1}}\fi}}}            
\providecommand{\constref}[2][]{\ensuremath{\constant{\textup{\ref{#2}{\ifx|#1|{}\else{,\ensuremath{#1}}\fi}}}}}
\providecommand{\constdef}[2][]{\label{#2}\ensuremath{\constant{\textup{\ref{#2}{\ifx|#1|{}\else{,\ensuremath{#1}}\fi}}}}}

\providecommand{\funkref}[3][]{\ensuremath{{#3}_{\textup{\ref{#2}{\ifx|#1|{}\else{,\ensuremath{#1}}\fi}}}}}
\providecommand{\funkdef}[3][]{\label{#2}\funkref[#1]{#2}{#3}}

\renewcommand{\div}[1][]{\nabla\ifx|#1|{}\else\kern-2pt_{#1}\fi\kern-2pt\inner}
\providecommand{\divof}[2][]{\div[#1]\ifx|#2|{}\else\qb{#2}\fi}
\providecommand{\grad}{}
\renewcommand{\grad}[1][]{\nabla\ifx|#1|\else_{#1}\fi}

\providecommand{\rot}[1][]{\nabla\ifx|#1|\else_{#1}\fi\outerp}

\providecommand{\rowdiv}[1][]{\D\ifx|#1|{}\else\kern-1pt_{#1}\kern-2pt\fi\cdot}
\providecommand{\rowdivof}[2][]{\rowdiv[#1]\ifx|#2|{}\else\qb{#2}\fi}

\providecommand{\inverse}[2][]{\powp[#1]{-1}{#2}}
\providecommand{\inverseqp}[1]{\inverse{\qp{#1}}}
\providecommand{\inverseof}[1]{\inverseqp{#1}}

\providecommand{\fracl}[3][]{\ifx|#1|\nicefrac{#2}{#3}\else{#2}#1/{#3}\fi}

\providecommand{\qpfracl}[3][]{\qp{\ifx|#1|\fracl{#2}{#3}\else{#2}#1/{#3}\fi}}
\providecommand{\qpfrac}[3][]{\qp{\ifx|#1|\frac{#2}{#3}\else{#2}#1/{#3}\fi}}
\providecommand{\absfracl}[3][]{\abs{\ifx|#1|\fracl{#2}{#3}\else{#2}#1/{#3}\fi}}
\providecommand{\absfrac}[3][]{\abs{\ifx|#1|\frac{#2}{#3}\else{#2}#1/{#3}\fi}}
\providecommand{\fraclff}[3][]{\ifx|#1|{#2}/{#3}\else{#2}#1/{#3}\fi}

\providecommand{\eye}[1][]{\vec{\mathrm I}\ifx|#1|{}\else_{#1}\fi}
\providecommand{\numeye}[1][]{\boldsymbol{\mathsf{I}}\ifx|#1|{}\else_{#1}\fi}
\providecommand{\Eye}[1]{
  \begin{bmatrix}
  \ifthenelse{#1>1}{
    \ifthenelse{#1>2}{
      \ifthenelse{#1>3}{
        \ifthenelse{#1>4}{
          1&\zeroentry&\dotso&\zeroentry
          \\
          \zeroentry&1&\dotso&\zeroentry
          \\
          \vdots&\vdots&\ddots&\vdots
          \\
          \zeroentry&\zeroentry&\dotso&1
        }{        
          1&\zeroentry&\zeroentry&\zeroentry
          \\
          \zeroentry&1&\zeroentry&\zeroentry
          \\
          \zeroentry&\zeroentry&1&\zeroentry
          \\
          \zeroentry&\zeroentry&\zeroentry&1
        }
      }{
        1&\zeroentry&\zeroentry
        \\
        \zeroentry&1&\zeroentry
        \\
        \zeroentry&\zeroentry&1
      }
    }{
      1&\zeroentry
      \\
      \zeroentry&1
    }
  }{
    1
  }
  \end{bmatrix}
}

\providecommand{\lebmeas}[1][]{\operatorname{l}^{#1}}     
\providecommand{\lebmeasof}[2][]{\ifx|#1|\left|#2\right|\else\lebmeas[#1]\qp{#2}\fi}         

\providecommand{\meshsize}[1][]{h\ifx|#1|\else_{#1}\fi}

\providecommand{\negpart}[1]{{#1}_-}           

\providecommand{\Oh} {\operatorname{O}}                   

\providecommand{\dash}[1][']{\ifthenelse{\equal{#1}{'}\OR\equal{#1}{''}}{#1}{^{(#1)}}}

\providecommand{\pdfrac}[2][]{\ensuremath{\frac{\partial\ifx|#1|\phantom{#2}\else{#1}\fi}{\partial{#2}}}} 
\providecommand{\pdfracpow}[3][]{\ensuremath{\frac{\partial^{#3}\ifx|#1|\phantom{#2}\else{#1}\fi}{\partial{#2}^{#3}}}} 

\providecommand{\pd}[2][]{\ensuremath{\partial_{#2}}{\ifx|#1|{}\else{\qb{#1}}\fi}} 

\providecommand{\pdt}[1][]{\pd[#1]t}                       
\providecommand{\pdx}[1][]{\pd[#1]x}                       
\providecommand{\pdz}[1][]{\pd[#1]z}                       

\renewcommand{\Im}{\operatorname{im}}                 
\renewcommand{\Re}{\operatorname{re}}                 
\providecommand{\imaginpart}[1][]{\Im{\ifx|#1|{}\else\qp{#1}\fi}} 
\providecommand{\realpart}[1][]{\Re{\ifx|#1|{}\else\qp{#1}\fi}} 

\providecommand{\transpose}{\intercal}

\providecommand{\transposed}{{}^\transpose}

\providecommand{\Transpose}[1]{\ensuremath{{#1}^{\transpose}}}

\providecommand{\Transposevec}[1]{\Transpose{\vec{#1}}}
\providecommand{\transposevec}[1]{\Transposevec{#1}}

\providecommand{\orthogonalto}[1][]{\ensuremath{\perp\ifx|#1|{}\else{_{#1}}\fi}}

\providecommand{\rowof}[1]{\ensuremath{\vecof{#1}}}

\providecommand{\colvec}[1]{\ensuremath{\vecof{#1}}}

\providecommand{\colvectwo}[2]{\ensuremath{\colvec{#1,#2}}}

\providecommand{\discolvec}[2][r]{\ensuremath{\begin{bmatrix}[#1]#2\end{bmatrix}}}

\providecommand{\discolvecthree}[4][r]{\ensuremath{\discolvec[#1]{#2\\#3\\#4}}}

\providecommand{\zeroentry}{\phantom0}

\providecommand{\smint}{\ensuremath{{\text{\textbf{/}}}\kern-.75em\smallint}}
\renewcommand{\smint}[1][]{\lower12.3pt\hbox{\begin{tikzpicture}\draw[line width=.75pt] (-3pt,-0.5)--(1pt,-0.5) node[pos=0.6]{$\int$};\path (-1.5pt,-24pt)node {\scriptsize $#1$};\end{tikzpicture}}}

\providecommand{\lap}{\ensuremath{\Delta}}
\providecommand{\lapin}[1][]{\lap\ifx|#1|\else_{#1}\fi}
\providecommand{\normalsymbol}{\operatorname{n}}
\providecommand{\normal}[1][]{\normalsymbol\ifx|#1|\else_{#1}\fi}
\providecommand{\normalto}[1]{\ensuremath{\normal[#1]}}
\providecommand{\normalder}[1][]{\ensuremath{\normal\ifx|#1|\else\qp{#1}\fi{\inner\grad}}}
\providecommand{\normalderto}[2][]{\ensuremath{\normalto{#2}\ifx|#1|\else\qp{#1}\fi{\inner\grad}}}

\providecommand{\tangentialsymbol}{\operatorname{t}}
\providecommand{\tangentialto}[2][]{\tangentialsymbol\ifx|#1|\else^{#1}\fi\ifx|#2|\else_{#2}\fi}

\providecommand{\union}[1]{\ensuremath{\bigcup}_{#1}}

\ifthenelse{\boolean{usesvjour}}{}{
  \let\vec\undefined
  \providecommand{\vec}[1]{\ensuremath{\boldsymbol{#1}}}
  \renewcommand{\vec}[1]{\ensuremath{\boldsymbol{#1}}}
}

\providecommand{\hatmat}[1]{\hat{\mat{#1}}}

\providecommand{\geomat}[1]{\vec{\MakeUppercase{#1}}}

\providecommand{\tildevec}[1]{\ensuremath{\widetilde{\vec{#1}}}}

\providecommand{\mat}[1]{\geomat{#1}} 

\providecommand{\Prob}[1][]{\ensuremath{\operatorname{Prob}\ifx|#1|{}\else_{#1}\fi}}

\providecommand{\pdf}[2][]{\ensuremath{\operatorname{pdf}_{#2\ifx|#1|{}\else{\conditionalto{#1}}\fi}}\xspace}

\providecommand{\expectation}{\ensuremath{\operatorname{E}}}
\providecommand{\EX}[1][]{\ensuremath{\expectation\ifx|#1|{}\else_{#1}\fi}}

\providecommand{\gausskernel}[3][x]{%
  \ensuremath{
    \exp\frac{-\if#20{#1}\else(#1-\mu)\fi^2}{%
      2\if#31{}\else\powp2{#3}\fi}%
  }%
}
\providecommand{\gaussdistribution}[3][x]{%
  \ensuremath{\frac1{\sqrt{2\pic}\if#31{}\else#3\fi}%
    \gausskernel[#1]{#2}{#3}
  }%
}%

\providecommand{\boundary}{\partial}

\providecommand{\SPD}{\operatorname{SPD}}
\providecommand{\spdmats}[2][F]{\SPD(\csname#1\endcsname{#2})}
 \providecommand{\Hspace}{\ensuremath{\operatorname H}\xspace}
 \providecommand{\Lebesgue}{\ensuremath{\operatorname L}\xspace}

 \providecommand{\Weaklyder}{\ensuremath{\operatorname W}\xspace}
 
 \providecommand{\dual}[1]{\ensuremath{{#1}'}}
 
 \providecommand{\dualspace}[2][]{\dual{\linspace{#2}\ifx|#1|\else{_{#1}}\fi}}
 \providecommand{\bidual}[1]{\ensuremath{{#1}''}}
 \providecommand{\bidualspace}[2][]{\bidual{\linspace{#2}\ifx|#1|\else{_{#1}}\fi}}

 \providecommand{\BV}[1]{\ensuremath{\operatorname{BV}}}
 
 \providecommand{\leb}[1]{\ensuremath{\Lebesgue_{#1}}}
 \providecommand{\lebloc}[1]{\ensuremath{{{\Lebesgue}^{\kern-.20em\lower .1ex\hbox{\tiny\textrm{\textup{loc}}}}_{#1}}}}
 \providecommand{\lebnorm}[3][]{\ensuremath{\Norm{#2}_{\leb{#3}\ifx|#1|{}\else(#1)\fi}}}

 \providecommand{\sob}[2]{\ensuremath{{\smash\Weaklyder}^{#1}_{#2}}}

 \providecommand{\sobh}[1]{\ensuremath{\Hspace^{#1}}}
 \providecommand{\vecsobh}[1]{\ensuremath{\vec\Hspace^{#1}}}
 \providecommand{\hdiv}[1][]{\vecsobh{\operatorname{div}}\ifx|#1|\else(#1)\fi}
 \providecommand{\hcurl}[1][]{\vecsobh{\operatorname{curl}}\ifx|#1|\else(#1)\fi}
 
 \providecommand{\sobhz}[2][]{\sobh{#2}_{0\ifx+#1+\else|#1\fi}}
 
 \providecommand{\Lip}[1][]{\ensuremath{\operatorname{Lip}}\ifx|#1|{}\else{\qp{#1}}\fi}

 \providecommand{\Symmatrices}[2][R]{\ensuremath{\operatorname{Sym}{(\csname#1\endcsname{#2})}}}
 
 \providecommand{\SAmatrices}[2][F]{\ensuremath{\operatorname{SA}{(\csname#1\endcsname{#2})}}}
 \providecommand{\mesh}[2][]{{\ensuremath{\mathcalbf{\MakeUppercase{#2}}\ifx|#1|\else_{#1}\fi}}}

\providecommand{\crouzeixraviart}[1][1]{\operatorname{CR}\ifx|#1|{}\else{^{#1}}\fi}

\providecommand{\linspace}[1]{\mathscript{\MakeUppercase{#1}}}

\providecommand{\clinopss}[2]{\clinopss{\linspace{#1}}{\linspace{#2}}}
\providecommand{\fepartition}[2][]{\mathscript{\MakeUppercase{#2}}\ifx|#1|{}\else_{#1}\fi}
\providecommand{\fespace}[2][]{\mathbb{\uppercase{#2}}\ifx|#1|{}\else_{#1}\fi}

\providecommand{\vespace}[1][]{\fespace v\ifx|#1|\else_{#1}\fi}

\providecommand{\fe}[2][]{\ensuremath{\uppercase{#2}\ifx|#1|\else_{#1}\fi}}

\providecommand{\vecfe}[2][]{\ensuremath{\vec{\uppercase{#2}}\ifx|#1|{}\else{_{#1}}\fi}}

\providecommand{\matfe}[2][]{\ensuremath{\mat{\uppercase{#2}}\ifx|#1|{}\else{_{#1}}\fi}}

\providecommand{\hatmatfe}[2][]{\ensuremath{\hatmat{\uppercase{#2}}\ifx|#1|{}\else{_{#1}}\fi}}

\providecommand{\Foreach}{\text{ for each }}

\RequirePackage{amscd}

\providecommand{\ideq}{\equiv}

\providecommand{\funk}[3]{\ensuremath{#1:#2\to#3}}

\providecommand{\dfunkmapsto}[6][]{\ensuremath{
    \begin{array}{rrcl}
      {#2}: & {#4} &  \to   & {#6}
      \\
            & {#3} &\mapsto & {#5\text{\ #1}}
    \end{array}\quad}}

\providecommand{\evalat}[3][]{\qb{#2}_{\ifx|#1|{}\else#1=\fi#3}}
\providecommand{\evaldiff}[4][]{\qb{#2}^{\ifx|#1|{}\else#1=\fi#3}_{\ifx|#1|{}\else#1=\fi#4}}

\providecommand{\aka}[1]{(also known as {#1})\xspace}

\providecommand{\akaindexemph}[2][]{\aka{\indexemph[#1]{#2}}}

\providecommand{\codename}[1]{\nolinkurl{#1}\xspace}

\ifthenelse{\boolean{useutopia}}{%
  
}{%
  
}

\providecommand{\indexen}[2][]{{\ifthenelse{\boolean{shownotes}}{\color b}{}#2\ifx|#1|\index{#2}\else\index{#1}\fi}}
\providecommand{\indexemph}[2][]{\emph{\indexen[#1]{#2}}}

\providecommand{\ListParameters}{}
\renewcommand{\ListParameters}%
{
	 \setlength{\topsep}{0pt}
	 \setlength{\leftmargin}{0pt}
         \setlength{\itemsep}{0pt}
	 \setlength{\parsep}{0pt}
	 \setlength{\parskip}{0pt}
         \setlength{\labelsep}{0pt}
	 \setlength{\itemindent}{0pt}
}

{%
  \begin{list}%
    {}%
    {\ListParameters%
    
}}%
{\end{list}}
\newcounter{tmpcounter}
\newcounter{LetterListItem}
\renewcommand{\theLetterListItem}{(\alph{LetterListItem})}
\newenvironment{LetterList}[1][0]%
{
	\begin{list}%
	{\theLetterListItem\ }%
	{\usecounter{LetterListItem}
	  \ListParameters
          \ifx|#1|{}\else\setcounter{LetterListItem}{#1}\fi
	}
}%
{\end{list}}
\newcounter{NumberListItem}
\renewcommand{\theNumberListItem}{\arabic{NumberListItem}}
{
	\begin{list}%
	{\theNumberListItem.\ }%
	{\usecounter{NumberListItem}%
	 \ListParameters
	}
}%
{\end{list}}
\newcounter{QuestionListItem}
\renewcommand{\theQuestionListItem}{\textbf{Question \arabic{QuestionListItem}}}
{
	\begin{list}%
	{\theQuestionListItem.\ }%
	{\usecounter{QuestionListItem}%
	 \ListParameters
	}
}%
{\end{list}}
\newcounter{RomanListItem}
\renewcommand{\theRomanListItem}{(\roman{RomanListItem})}
{
	\begin{list}%
	{\theRomanListItem\ }%
	{\usecounter{RomanListItem}
	 \ListParameters
	}
}%
{\end{list}}
\newcounter{StepsItem}
{
	\begin{list}%
	{Step \theStepsItem.\ }%
	{\usecounter{StepsItem}%
	 \ListParameters
	}
}%
{\end{list}}
\newcounter{CasesListItem}
\renewcommand{\theCasesListItem}{\Alph{CasesListItem}}
{
	\begin{list}%
	{\emph{Case \theCasesListItem.}\ }%
	{\usecounter{CasesListItem}%
	 \ListParameters
	}
}%
{\end{list}}
\newcounter{QAListItem}
\renewcommand{\theQAListItem}{Q\arabic{QAListItem}:}
{
	\begin{list}%
	{\theQAListItem}%
	{\usecounter{QAListItem}
	 \ListParameters
	}
}%
{\end{list}}

\ifthenelse{\boolean{isthesis}}{%
  \setcounter{secnumdepth}{1}%
}{
  \setcounter{secnumdepth}{2} 
}
\providecommand{\ListParameters}{}
\renewcommand{\ListParameters}
{
	 \setlength{\topsep}{0em}
	 \setlength{\leftmargin}{0em}
         \setlength{\itemsep}{0ex}
	 \setlength{\parsep}{.5ex}
	 \setlength{\itemindent}{\labelsep}
	 \addtolength{\itemindent}{\labelwidth}
}

  \providecommand{\ObsName}{Remark}
  \providecommand{\DefName}{Definition}
  \providecommand{\ExaName}{Example}
  \providecommand{\TheName}{Theorem}
  \providecommand{\LemName}{Lemma}
  \providecommand{\ProName}{Proposition}
  \providecommand{\CorName}{Corollary}
  \providecommand{\PbmName}{Problem}
  \providecommand{\AlgName}{Algorithm}
  \providecommand{\ExeName}{Exercise}
  \providecommand{\SolName}{Solution}

\ifthenelse{\boolean{isthesis}}{%
  
}{%
  
}
\makeatletter
\newcommand{\oltikzgetxy}[3]{%
  \tikz@scan@one@point\pgfutil@firstofone#1\relax
  \edef#2{\the\pgf@x}%
  \edef#3{\the\pgf@y}%
}
\makeatother
\providecommand{\pdfformat}[1]{
   \provideboolean{pdfoutput}
   \setboolean{pdfoutput}{#1}
  \ifthenelse{\boolean{pdfoutput}}{
    \typeout{using pdf}
\usepackage{pdfsync}
    \providecommand{\graphext}{pdf}
    \renewcommand{\graphext}{pdf}
    \providecommand{\graphextex}{pdf_t}
    \renewcommand{\graphextex}{pdf_t}
  }{
    \typeout{using eps}
    \RequirePackage[dvips]{graphicx,xcolor}
    \providecommand{\graphext}{eps}
    \renewcommand{\graphext}{eps}
    \providecommand{\graphextex}{eps_t}
    \renewcommand{\graphextex}{eps_t}
  }
  \RequirePackage{epsfig}
  \RequirePackage{tikz}
  \RequirePackage{rotating}
  
  \RequirePackage{graphicx}
  \RequirePackage{xcolor}
  \provideboolean{darkcolortheme}
  \definecolor{SussexFlint}{rgb}{.00,.19,.21}
  \definecolor{SussexGrey}{rgb}{.51,.58,.49}
  \definecolor{SussexOrange}{rgb}{.94,.29,.00}
  \definecolor{SussexYellow}{rgb}{1.00,.73,.00}
  \definecolor{SussexRed}{rgb}{.94,.01,.49}
  \definecolor{SussexPurple}{rgb}{.48,.06,.44}
  \definecolor{SussexGreen}{rgb}{.00,.58,.46}
  \definecolor{OmarGreen}{rgb}{.00,.68,.36}
  \definecolor{SussexBlue}{rgb}{.00,.58,.65}
  \definecolor{OmarBlue}{rgb}{.00,.38,.65}
  \colorlet{a}{OmarBlue}
  \colorlet{b}{SussexOrange}
  \colorlet{c}{SussexGreen}
  \colorlet{d}{SussexPurple}
  \colorlet{e}{SussexRed}
  \colorlet{f}{SussexYellow}
  \colorlet{g}{white}
  \colorlet{h}{SussexGrey}
  \colorlet{i}{black}
  \colorlet{j}{SussexFlint}
  \colorlet{colora}{a}
  \colorlet{colorb}{b}
  \colorlet{colorc}{c}
  \colorlet{colord}{d}
  \colorlet{colore}{e}
  \colorlet{colorf}{f}
  \colorlet{colorg}{g}
  \colorlet{colorh}{h}
  \colorlet{colori}{i}
  \colorlet{colorj}{j}
  \newcommand{\mausDarkColorTheme}{
    \colorlet{a}{SussexYellow!50!yellow}
    \colorlet{b}{SussexBlue}
    \colorlet{c}{SussexRed!50!red}
    \colorlet{d}{SussexOrange!50!yellow}
    \colorlet{e}{SussexGreen!50!green}
    \colorlet{f}{SussexPurple!50!magenta}
    \colorlet{g}{black}
    \colorlet{h}{SussexFlint!50!black}
    \colorlet{i}{white}
    \colorlet{j}{SussexGrey}
  }
  \ifthenelse{\boolean{darkcolortheme}}{\mausDarkColorTheme}{}
}
\providecommand{\solution}{\textbf{\SolName.}\xspace}

 \newcounter{phantombox}[enumi]
 \provideboolean{showphantoms}
 \renewcommand{\thephantombox}{\Alph{phantombox}}
 \newcommand{\phantombox}[1]{\stepcounter{phantombox}%
   \ensuremath{\boxed{%
       {\ifthenelse{\boolean{showphantoms}}{#1}{\phantom{#1}}}%
       {\texttt{\tiny\ \colorbox{i!50}{\color g\thephantombox}}
       }%
     }%
   }%
 }

 \provideboolean{hidesolution}
 \newcommand{\consolution}[2][]{
   \ifthenelse{\boolean{hidesolution}}{#1\setboolean{showphantoms}{false}}{%
     {\setboolean{showphantoms}{true}\color{i!50}\par \small {\solution}\ #2\par\ \\[5pt]}}
 }
 \provideboolean{showmarks}
 \providecommand{\showmarks}[1]{
   \ifthenelse{
     \boolean{showmarks}}{
     \marginpar{
       \tiny [$#1$ mark\ifthenelse{\equal{#1}1}{\phantom{s}}s]}%
   }{}}

 \newcommand{\condibreak}{\ifthenelse{\boolean{hidesolution}}{\newpage}{}}
 
 \providecommand{\qeyword}[1]{\index{#1}\ifthenelse{\boolean{shownotes}}{\texttt{\color{e}[#1]}}{}}
 \providecommand{\sourcecite}[2][]{\index{#1}\ifthenelse{\boolean{shownotes}}{\texttt{\color{d}[source: \cite[#1]{#2}]}}{}}
 
 \RequirePackage{hyphenat}
\usepackage{xspace}
\usepackage{ifthen}
\usepackage{iftex}
\usepackage[utf8]{inputenc}
\usepackage{siunitx}
\usepackage{url}
\usepackage{amsmath}
\usepackage{amssymb}
\usepackage{latexsym}
\usepackage{enumerate}
\usepackage{graphicx}
\usepackage{array}
\usepackage{parskip}
\usepackage{empheq}
\usepackage{mathtools}
\usepackage{xfrac}
\usepackage[normalem]{ulem}
\usepackage{textcomp}
\usepackage{framed}
\usepackage{multicol}
\provideboolean{shownotes}
\usepackage[justification=justified]{caption}
\usepackage{subfigure}
\usepackage{float}
\usepackage{totcount}
\usepackage{booktabs}
\usepackage{marginnote}
\usepackage{titleps}
\usepackage{enumitem}

\usepackage{tikz}
\usetikzlibrary{arrows, automata, calc, patterns, positioning, shapes}
\usetikzlibrary{decorations.text}
\usepackage{etoolbox}
\usepackage{xcolor, colortbl}
\usepackage{pgfplots}
\usepgfplotslibrary{dateplot, patchplots, statistics}

\providecommand{\indexen}[1]{#1\index{#1}}
\providecommand{\indexemph}[1]{\emph{\indexen{#1}}}
\providecommand{\funk}[3]{\ensuremath{#1:#2\to#3}}
\providecommand{\reals}{\rR}
\providecommand{\rR}{\ensuremath{\mathbb R}\xspace}
\renewcommand{\d}{\ensuremath{\operatorname{d}\!\!}}
\providecommand{\xt}{(t,x)}

\providecommand{\vect}[1]{\boldsymbol{#1}}
\providecommand{\tildevec}[1]{\ensuremath{\widetilde{\vect{#1}}}}

\providecommand{\gravity}{\ensuremath{\operatorname g}}

\providecommand{\BigO}[1]{\ensuremath{\operatorname{O}\!\left(#1\right)}}
\providecommand{\aka}[1]{(also known as {#1})\xspace}
\providecommand{\varep}{\varepsilon}
\providecommand{\Reynolds}{\text{Rey}}
\providecommand{\Frhor}{\text{Fro}}

\providecommand{\funkref}[3][]{\ensuremath{{#3}_{\textup{\ref{#2}{\ifx|#1|{}\else{,\ensuremath{#1}}\fi}}}}}
\providecommand{\funkdef}[3][]{\label{#2}\ensuremath{{#3}_{\textup{\ref{#2}{\ifx|#1|{}\else{,\ensuremath{#1}}\fi}}}}}

\providecommand{\uu}{\vect{u}}
\providecommand{\vv}{\vect{v}}
\providecommand{\W}{\Omega}
\providecommand{\ue}{u_{\varep}}
\providecommand{\ve}{v_{\varep}}
\providecommand{\pe}{p_{\varep}}

\providecommand{\average}[1]{\left<#1\right>}

\providecommand{\condicolor}[1][a]{\ifthenelse{\boolean{shownotes}}{\color{#1}}{}}
\numberwithin{equation}{subsection}

\providecommand{\ListParameters}{}
\renewcommand{\ListParameters}%
{
	 \setlength{\topsep}{0pt}
	 \setlength{\leftmargin}{0pt}
         \setlength{\itemsep}{0pt}
	 \setlength{\parsep}{0pt}
	 \setlength{\parskip}{0pt}
         \setlength{\labelsep}{0pt}
	 \setlength{\itemindent}{0pt}
}

\providecommand{\qp}[1]{\ensuremath{\left({#1}\right)}}
\providecommand{\qb}[1]{\ensuremath{\!\left[{#1}\right]}}
\providecommand{\qc}[1]{\ensuremath{\left\{{#1}\right\}}}
\providecommand{\pd}[2][]{\ensuremath{\partial_{#2}}{\ifx|#1|{}\else{\qb{#1}}\fi}\xspace}
\providecommand{\pdt}[1][]{\pd[#1]t}
\providecommand{\pdx}[1][]{\pd[#1]x}

\providecommand{\pdz}[1][]{\pd[#1]z}
\providecommand{\fracl}[3][]{{#2}#1/{#3}}
\providecommand{\setof}[1]{{\qc{#1}}}
\providecommand{\inverse}[1]{\ensuremath{{#1}^{-1}}}
\providecommand{\inverseqp}[1]{\inverse{\qp{#1}}}
\providecommand{\inverseof}[1]{\inverseqp{#1}}
\providecommand{\const}{\operatorname{constant}}
\providecommand{\constant}[1]{\ensuremath{C_{#1}}}
\providecommand{\constext}[2][]{\constant{\text{#2}{\ifx|#1|{}\else{,\ensuremath{#1}}\fi}}}
\providecommand{\discolvec}[2][r]{\ensuremath{\begin{bmatrix}[#1]#2\end{bmatrix}}}

\providecommand{\discolvecthree}[4][r]{\ensuremath{\discolvec[#1]{#2\\#3\\#4}}}

\providecommand{\one}{\ensuremath{\mathbb 1}\xspace}
\providecommand{\charfun}[1]{\one_{#1}}
\providecommand{\iverson}[1]{\one\qb{#1}}
\providecommand{\abs}[2][]{\ensuremath{\ifx|#1|{\left|#2\right|}\else{\csname#1\endcsname|{#2}\csname#1\endcsname|}\fi}}
\providecommand{\tand}{\ensuremath{\text{ and }}}

\providecommand{\diver}[1][]{\text{div} \ \ifx|#1|{}\else\kern-3pt_{#1}\fi\kern-3pt}
\providecommand{\divof}[2][]{\diver[#1]\ifx|#2|{}\else\qb{#2}\fi}

\providecommand{\transpose}{\intercal}
\providecommand{\transposed}{{}^\transpose}
\providecommand{\geomat}[1]{\vect{#1}}
\providecommand{\mat}[1]{\geomat{#1}}
\providecommand{\vecof}[1]{\qp{#1}}
\providecommand{\colvec}[1]{\ensuremath{\vecof{#1}}}
\providecommand{\colvectwo}[2]{\ensuremath{\colvec{#1,#2}}}

\providecommand{\ideq}{\equiv}
 \providecommand{\opinter}[2]{\ensuremath{\left(#1,#2\right)}\xspace}
 \providecommand{\clinter}[2]{\ensuremath{\left[#1,#2\right]}\xspace}
 \providecommand{\dfunkmapsto}[6][]{\ensuremath{
    \begin{array}{rrcl}
      {#2}: & {#4} &  \to   & {#6}
      \\
            & {#3} &\mapsto & {#5\text{\ #1}}
    \end{array}\quad}}
\providecommand{\ensemble}[2]{\ensuremath{\left\{ #1:\;#2 \right\}}}
\providecommand{\union}[1]{\ensuremath{\bigcup}_{#1}}
\renewcommand{\vec}[1]{\ensuremath{\boldsymbol{#1}}}
\providecommand{\Transpose}[1]{\ensuremath{{#1}^{\transpose}}}
\providecommand{\Transposevec}[1]{\Transpose{\vec{#1}}}
\providecommand{\transposevec}[1]{\Transposevec{#1}}
 \providecommand{\getrow}[2]{\irow{\entry{#1}}{#2}}

\providecommand{\fricI}[1][I]{k_-(#1)}
\providecommand{\fricR}[1][R]{k_+(#1)}
\providecommand{\gravity}{\ensuremath{\operatorname g}}
\newcommand{\LHS}[1]{\text{LHS} \eqref{#1}}
\newcommand{\RHS}[1]{\text{RHS} \eqref{#1}}
\newcommand{\CFL}{\ensuremath{\operatorname{CFL}}}
\renewcommand{\div}[1][]{\text{div} \ \ifx|#1|{}\else\kern-3pt_{#1}\fi\kern-3pt}
\renewcommand{\one}{\ensuremath{\mathbb 1}\xspace}

\providecommand{\qed}{\vrule height 5pt depth 0pt width 3pt}
\providecommand{\qqed}{{\raggedright{\ \hfill\qed}}}
\newenvironment{Proofof}[1]{\subsection{Proof of #1}}{\qqed\par}
\setboolean{showchanges}{false}
\setboolean{shownotes}{false}
\usepackage{natbib}
\usepackage{backref}
\usepackage{lipsum}
\makeindex
\usepackage{dsfont,bm,stmaryrd}

\providecommand{\ThesisSimplify}[3]{
\ifthenelse{#1}{\begin{framed}\centering \textbf{{\color{red} #2 suppressed to speed up compiling.}}\end{framed}}{#3}}
\newboolean{simplify}
\newboolean{nosimplify}
\setboolean{simplify}{false}
\setboolean{nosimplify}{false}

\usepackage{url}
\usepackage{amsthm}
\newtheorem{theorem}{Theorem}
\newtheorem{remark}{Remark}
\newtheorem{proposition}{Proposition}
\newtheorem{property}{Property}
\newtheorem{corollary}{Corollary}
\newenvironment{The}[1][]{\subsection{\TheName\xspace{\ifx&#1&{}\else{ (#1)}\fi}}\slshape}{\upshape}
\newenvironment{Pro}[1][]{\subsection{\ProName\xspace{\ifx&#1&{}\else{ (#1)}\fi}}\slshape}{\upshape}

\newenvironment{Cor}[1][]{\subsection{\CorName\xspace{\ifx&#1&{}\else{ (#1)}\fi}}\slshape}{\upshape}

\newenvironment{Pty}[1][]{\subsection{Property\xspace{\ifx&#1&{}\else{ (#1)}\fi}}\slshape}{\upshape}

\newenvironment{Obs}[1][]{\subsection{\ObsName\xspace{\ifx&#1&{}\else{ (#1)}\fi}}}{}
\newenvironment{Exa}[1][]{\subsection{\ExaName\xspace{\ifx&#1&{}\else{ (#1)}\fi}}}{}

\numberwithin{equation}{subsection}
\hypersetup{
    unicode=true,		
    pdftoolbar=true,        	
    pdfmenubar=true,        	
    pdffitwindow=true,      	
    pdftitle={A Saint-Venant model for overland flows withOA rain},    
    pdfauthor={Mehmet Ersoy and Omar Lakkis and Philip Townsend},     
    pdfsubject={mathematics, analysis, partial differential equation, numerical, computational},   
    pdfcreator={pdflatex},   	
    pdfproducer={}, 		
    pdfkeywords={numerical analysis, computational method, conservation law, finite volume, kinetic equation, Saint Venant, shallow water}, 			
    pdfnewwindow=true,
    colorlinks=false,
    linkcolor=blue,  
    citecolor=blue,        	
    filecolor=magenta,      	
    urlcolor=cyan           	
}
\title[A Saint-Venant model with rain]{A Saint-Venant model for overland flows with precipitation and recharge}
\author{Mehmet Ersoy}
\address{
  Mehmet Ersoy
  \newline
  Université de Toulon
  \\
  IMATH EA 2134
  \\
  La Garde (France)
  \\
  FR-83957%
}
\email{\linkedemail{Mehmet.Ersoy@univ-tln.fr}}
\urladdr{\linkedurl{http://ersoy.univ-tln.fr/}}
\author{Omar Lakkis}
\address{ Omar Lakkis
  \newline
  Department of Mathematics
  \\
  University of Sussex
  \\
  Brighton (England, United Kingdom)
  \\
  GB-BN1 9QH%
}
\curraddr{}
\email{\linkedemail{lakkis.o.maths@gmail.com}}
\author{Philip Townsend}
\address{
  Philip Townsend 
  \newline
  Chalmers University of Technology
  \\
  Chalmersplatsen 4
  \\
  412 96 Gothenburg (Sweden)
}
\curraddr{}
\email{\linkedemail{phil.townsend123@gmail.com}}
\renewcommand{\d}{\ensuremath{\operatorname d}}
\begin{document}
\date{}
\maketitle
\begin{abstract} 
We propose a one-dimensional Saint-Venant (open channel) model
overland flows including a water input--output source term modelling
recharge via rainfall and infiltration (or exfiltration). We derive
the model via asymptotic reduction from the two-dimensional
Navier--Stokes equations under the shallow water assumption, with
boundary conditions including recharge via ground infiltration and
runoff. This new model recovers existing models as sepcial cases, and
adds more scope by adding a water-mixing friction terms that depends
on the rate of water recharge. We prospose a novel entropy function
and its flux, that are useful in validating the model's conservation
or dissipation properties. Based on this entropy function we propose a
finite volume scheme extending a class of kinetic schemes and provide
numerical comparisons with respect to the newly introduced mixing
friction coefficient. We also provide a comparison with experimental
data.
\end{abstract}
\clearpage
\clearpage
\section{Introduction}
\label{SectionIntro}
\renewcommand{\theequation}{\arabic{section}.\arabic{equation}}
\changes{In quantifying the dynamics of a watercourse, the most
  important component of the hydrologic recharge and loss is the the
  \emph{precipitation} and \emph{infiltration} processes,
  respectively. This is particularly important today in understanding
  and forecasting the impact of climate variability on the human and
  natural environment.} Modelling these processes and predicting the
  motion of water is a difficult task to which substantial effort has
  been devoted
  \citep{
    GraceEagleson:66:article:Modeling,
    WoolhiserLiggett:67:article:Unsteady,
    ZhangCundy:89:article:Modeling,
    EstevesFaucherGalleVauclin:00:article:Overland,
    WeillMouchePatin:09:article:A-generalized,
    RousseauCerdanErnLe-Maitre:12:article:Study
  }.
  
One of the most widely used models to describe the overland motion of
watercourses is the classical \emph{one-dimensional} \indexemph{Saint-Venant
  system} (also known as the \indexemph{open channel} or \indexemph{shallow water equations})
developed by Adhémar Jean Claude Barré de Saint-Venant in the 19th
Century as a reduction of the Navier--Stokes equation under certain
assumptions on the horizontal and vertical scales
\citep{Saint-Venant:71:article:Theorie}. For the specific problem of
modelling flooding caused by \indexen{precipitation},
\changes{the inclusion of
a source term corresponding to the recharge or infiltration
in the Saint-Venant system turns it from a conservation
law into a balance law.}
Existing approaches modelling surface flows under the effect of
rainfall or runoff are provided, for example, by
\citet{Sochala:08:phdthesis:Numerical} and
\citet{DelestreCordierDarbouxJames:12:article:A-limitation},
who model this phenomenon using the system
\begin{equation}
  \label{eq:rainy-Saint-Venant-according-to-Delestre}
  \begin{aligned}
    \pdt  h + \pdx\qb{hu} 
    &= S\\
    \pdt  \qb{hu} 
    + \pdx\qb{hu^2 + \frac{\gravity h^2}{2}}
    &= 
    -\gravity  h\pdx  Z - k_0(u)u
  \end{aligned}
\end{equation}
where the unknowns~\index{$h$}$h(t,x)$ and~\index{$u$}$u(t,x)$ model,
respectively, the \emph{height} and \emph{velocity of the water
  column} at space-time point $(t,x)$, $\gravity$ the gravitational
acceleration (considered a constant $\approx 9.81 m/s^{2}$), $Z(x)$
the topography of the channel bed with slope $\pdx Z(x)$, and $k_0$ an
empirical fluid-wall friction.  The \indexemph{source term}
\indexen{$S$} quantifies the amount of water that is added to ($S >
0$) or subtracted from ($S < 0$) from the flow, which in practice may
occur through a variety of mechanisms, (e.g. \indexen{direct
  rainfall}, \indexen{lateral flow}, \indexen{run-off},
\indexen{smaller tributaries}).  Among early works including the effect of rainfall
(or lateral inflow) on surface flows that relate to this
research, we note \citet{WoolhiserLiggett:67:article:Unsteady},
\citet{Wenzel:70:techreport:The-Effect}, and
\citet{ZhangCundy:89:article:Modeling}.


Our goal in this paper is to derive a model akin to
\eqref{eq:rainy-Saint-Venant-according-to-Delestre} via vertical
averaging under the shallow water assumption, starting from the
Navier--Stokes equations with a permeable Navier boundary condition to
account for the infiltration and a kinematic boundary condition to
consider the precipitation. The obtained averaged model extends this
system in a unique manner through an additional momentum source term
of the form
\begin{equation}
  \begin{aligned}
    Su - (\fricR + \fricI)u\quad\text{ with } S := R - I,
  \end{aligned}
\end{equation}
where $R \geq 0$ denotes the \emph{recharge rate} on the free surface
(accounting for both rain and runoff effects) and $I$ denotes the
\emph{infiltration rate} from the water to the ground (when $I > 0$)
or the ground to the water (when $I < 0$), i.e., \emph{seepage},
sometimes called \emph{exfiltration}. The terms $\fricR$ and $\fricI$,
which will be discussed in detail in \S\ref{sec:geometric-setup},
model the friction caused by \indexemph{recharge}, i.e., the addition
of water (assumed to have zero horizontal velocity), which attaches to
and is advected by the flow.  We will see below that these friction
terms are necessary to avoid paradoxical outcomes such as perpetual
motion, and, for simplicity, we will assume in this paper the most
basic constitutive relations for this friction: linear in $R$ for
$\fricR$ and piecewise linear in $I$ for $\fricI$. We note, however,
that these terms could be generalised by having two separate friction
coefficients or by replacing the linearity with more precise
constitutive relations.

We outline the rest of the article as follows: in
\S\ref{SectionModelPhy}, we present the geometric setup of the system
and the adjusted boundary conditions (including precipitation,
infiltration, and the corresponding friction terms) of the typical
Navier--Stokes equations. In \S\ref{SectionAveragedModel}, we derive
the consequent Saint-Venant system through a first order approximation
and discuss several theoretical results and corollaries that can be
derived for the system in \S\ref{sec:entropy-condition}. In
\S\ref{SectionKinSch}, we adapt to our model the finite volume kinetic scheme
considered in \cite{AudusseBristeauPerthame:00:techreport:Kinetic} and
\cite{PerthameSimeoni:01:article:A-kinetic}, and finally
present numerical experiments of the resulting code to
demonstrate the application of the model in \S\ref{SectionNum}. A C
and C++ implementation of this code, written by Matthieu Besson, Omar
Lakkis and Philip Townsend, is freely available on request (an older
version is given by \cite{BessonLakkis:13:url:Finite}).

\renewcommand{\theequation}{\arabic{section}.\arabic{subsection}.\arabic{equation}}
\section{Navier--Stokes equations with infiltration and recharge}
\label{SectionModelPhy}
Our aim is to construct a mathematical model for overland flows that is consistent with the physical phenomena that can affect the motion of such water. To this purpose, we propose a model reduction of the two-dimensional Navier--Stokes equations leading to an extension of the standard Saint-Venant system. By considering suitably chosen boundary conditions, we take into account the addition and removal of water, either by rainfall (e.g. from runoff onto the top of the watercourse) or by groundwater infiltration or exfiltration processes (e.g via a porous soil).

We start in~\S\ref{sec:geometric-setup} by reviewing the Navier--Stokes equations in the special geometric setting, describing the physics with a \indexemph{wet boundary} on the bottom of the water course and a \indexemph{free surface} on the top. We then introduce the boundary conditions for each surface in~\S\ref{sec:wet-boundary} and~\S\ref{sec:free-boundary}, respectively.

\subsection{Geometric set-up and the two-dimensional Navier--Stokes equations}
\label{sec:geometric-setup}
With numerical and practical applications in mind, we assume an arbitrary final time $T>0$. With reference to Fig.~\ref{fig:river-and-river-bed}, we consider an incompressible fluid moving in the space-time box
\begin{equation}
\begin{gathered}
  [0,T] \times \R2 \text{ with typical point denoted } (t,x,z)\\
\end{gathered}
.
\end{equation}
The \indexemph{absolute height} of the surface of the watercourse and the \indexemph{topography} of the channel bed are modelled, respectively, by the functions
\begin{equation}
  \dfunkmapsto[,]H{(t,x)}{\clinter0T\times\reals}{H(t,x)}{\reals}
  \quad
  \dfunkmapsto[,]Zx{\reals}{Z(x)}\reals 
\end{equation}
whose values measure with respect to a \indexemph{reference horizontal height} $0$. We define the \indexemph{local height} of the water by
\begin{equation}
  \label{eq:localwaterelevation}
  h(t, x) := H(t, x) - Z(x).
\end{equation}
\emph{We assume that there exists $c\in\reals$ such that $\int_\reals
  Z(x)-c\d x$ is finite.}  This assumption which is not restrictive
from a modelling point of view will be necessary for entropy balance
considerations.  Furthermore the vertical reference axis can be chosen
so that $c=0$.
The \indexemph{wet region} is defined as the area in which the fluid
resides at each time $t\in\clinter0T$
\begin{equation}
  \W(t) := \ensemble{(x, z) \in \reals^2}{Z(x) < z < H(t, x)},
\end{equation}
with its global counterpart
\begin{equation}
  \W := \union{0\leq t\leq T} \W(t)
  .
\end{equation}

\begin{figure}[ht]
  \centering    
  \scalebox{.8}{\begin{tikzpicture}
      \providecommand{\tikzinterface}{\x*\x*1/128-(\x-5)/8+1}
      \begin{scope}
        \clip (0,3.5) -- plot[smooth, domain=0:10]({\x}, {\tikzinterface}) -- plot[smooth, domain=10:0]({\x}, {0.2*cos(deg(\x))+4.5});
        \fill [cyan,opacity=0.3] (-0.5,-0.5) rectangle (10,5);
      \end{scope}
      
      \begin{scope}
        \clip (0,0) -- plot[smooth, domain=0:10]({\x}, {\tikzinterface}) -- (10,4) -- (10,0);
        \fill [brown,opacity=0.3] (-0.5,-0.5) rectangle (10,4);
      \end{scope}
      
      \draw [cyan, thick] plot[smooth, domain=0:10]({\x}, {0.2*cos(deg(\x))+4.5});
      \draw [domain=0:10,smooth,variable=\x,brown,thick] plot({\x}, {\tikzinterface});
      \draw [black, -latex] (-0.5,0) -- (10.5,0);
      \draw [black, -latex] (0,0) -- (0,5.5);
      
      \draw [black, |-latex] (-0.25,0)--(-0.25,1.195) node[pos=0.5,sloped,above] {$Z(x)$};
      \draw [black, |-latex] (-0.25,1.195)--(-0.25,4.55) node [pos=0.5,sloped,above] {$h\xt$};
      \draw [black, -latex] (1.5,2.75) -- node [above] {$q\xt$} (3,2.75);
      \draw [cyan] (7.75,2.75) node [above] {wet region};
      \draw [brown] (7.75,0.5) node [above] {ground};
      
      \draw [cyan, ultra thick, fill=cyan!50] (4.9,1.195) rectangle (5.1,4.55);
      \draw [cyan] (4.9,2.8725) node [above, rotate=90] {water column};
      
      \draw [black, dashed] (5,1.195) -- (0,1.195);
      \draw [black, dashed] (5,4.55) -- (0,4.55);
      \draw [fill=black] (5,1.195) circle (0.05cm);
      \draw [fill=black] (5,4.55) circle (0.05cm);
      
      \draw [purple, thick, -latex] (5,4.55) -- node [below right] {$R\xt$} (5,4);
      \draw [black, -latex] (5,4.55) -- ++(10:1) node [pos=1,sloped,above] {$\tangentialto{\Omega}$};
      \draw [black, -latex] (5,4.55) -- ++(100:1) node [left] {$\normalto{\Omega}$};
      
      \draw [purple, thick, -latex] (5,1.195) -- node [right] {$I\xt$} (5,0.645);
      \draw [black, -latex] (5,1.195) -- ++(177:1) node [above] {$\tangentialto\Omega$};
      \draw [black, -latex] (5,1.195) -- ++(267:1) node [right] {$\normalto{\Omega}$};
      
      \draw (0,0) node [below] {$0$};
      \draw [fill=black] (5,0) circle (0.05cm) node [below] {$x$};
      \draw (10,0) node [below] {$L$};
      
      \draw [cyan] (10,4.35) node [right] {free surface $\mathcal{F}$};
      \draw [brown] (10,1.781) node [right] {wet boundary $\mathcal{B}$};
  \end{tikzpicture}}
  \caption{Diagram of a river and river bed depicting the variables of interest}
  \label{fig:river-and-river-bed}
\end{figure}

As we can see in Figure~\ref{fig:river-and-river-bed}, the wet region
has two boundaries; the first is the \indexemph{wet boundary} between the wet region and the ground, denoted by\index{\(\mathcal B\)}
\begin{equation}
  \mathcal{B} = \{(x,Z) : x \in\reals\},
\end{equation}
and the second is the \indexemph{free surface} between the wet region and
the surrounding air, denoted by\index{\(\mathcal F\)}
\begin{equation}
  \mathcal{F} = \{(t,x,H) : t > 0, x \in\reals\}.
\end{equation}

We assume that the viscous flow $\uu$ satisfies, on the space-time domain $\W$, the two-dimensional incompressible Navier--Stokes equation
\begin{equation}
  \label{eq:incompressible-2d-Navier--Stokes}
  \begin{gathered}
    \divof{\rho_0 \uu\transposed} = 0,
    \\
    \pdt\qb{\rho_0 \uu} + \divof{\rho_0 \uu\otimes \uu} 
    - \div{\mat\sigma\qb{\uu}} - \rho_0\vec F = 0
  \end{gathered}
\end{equation}
where $\uu=\colvectwo{u}{\vv}$ is the velocity field, $\rho_0$ is the
density of the fluid (taken to be constant since the fluid is incompressible),
$\vec{F}=\colvectwo{0}{-\gravity}$ is the external force of gravity
with constant $\gravity$, and $\mat\sigma\qb{\uu}$ is the total stress
tensor whose matrix given by
\begin{equation}
  \label{NSTotalStressTensor}
  \mat\sigma\qb{\uu} 
  :=
  \begin{bmatrix}
    - p+ 2 \mu \pdx  u             
    &
    \mu \qp{\pd zu + \pd x\vv}
    \\
    \mu \qp{\pd zu + \pd x\vv}  
    &
    - p + 2 \mu \pd{z} \vv
  \end{bmatrix}
\end{equation}
where $p$ is the pressure and $\mu > 0$ the dynamic viscosity. The
(algebraic) tensor product of two vectors $\vec{a}\otimes\vec{b}$ is
defined as $\vec a\transposevec b$ (all vectors are displayed as
columns) and the $\div$ of a covector/tensor is taken as the row-wise
divergence of the associated matrix; in coordinates this means
\begin{equation}
  \getrow{\div\mat\alpha}i = \sum_{j=x,z}\pd j\matentry\alpha ij
  \text{ for }i=x,z
\end{equation}

To work with the wet region, we introduce its \indexemph{indicator function}
\begin{equation}
  \Phi(t,x,z)
  := \charfun{\W(t)}(x, z)
  = \iverson{Z(x) \leq z \leq H(t, x)} \text{ for all } t, x, z \in \reals.
\end{equation}
The function $\Phi$ is advected by the flow so its material derivative, with respect to the flow $\uu$, must therefore be zero. Moreover, thanks to the incompressibility condition, $\Phi$ satisfies the following \indexemph{indicator transport equation}
\begin{equation}
  \label{eqn:indicator-transport}
  \pdt \Phi + \pdx\qb{\Phi u} + \pdz\qb{\Phi \vv} = 0 \text{ on } \W.
\end{equation}


\subsection{The wet boundary}
\label{sec:wet-boundary}
Crucial to our model derivation is the particular situation on the
\emph{wet boundary}, where the effect of \indexemph{infiltration}
plays a central role. Given a set $G \in \reals^2$ and a point $\vec
x\in\boundary G$, we denote by $\tangentialto G(\vec x)$ the unique
normalized tangential vector and by $\normal[G](\vec x)$ its outward
boundary normal (see Fig.~\ref{fig:river-and-river-bed} for $G = \W$).

On the wet boundary, the topography is assumed to be rough and hence
produces friction, which we take into account by considering the
following \emph{Navier boundary condition}:
\begin{equation}
  \label{eq:Navier-boundary-condition}
  \qp{\mat\sigma\qb{\uu}\normal[\W]} 
  \inner 
  \tangentialto \W 
  = 
  -\rho_0 \qp{k(\uu) - \fricI}\uu
  \inner 
  \tangentialto \W \quad \text{on } \mathcal{B}
\end{equation}
We will leave defining $\fricI$ for the moment and note that the
scalar function $k(\uu)$ models a general \emph{kinematic friction
  law} on the channel bed:
\begin{equation}
  \label{Friction}
  k(\vec\xi)
  := 
  (\constext{lam} + \constext{tur}\norm{\vec\xi})
  ,
  \quad \text{for all }\vec\xi\in\reals^2
\end{equation}
where the friction coefficients $\constext{lam}$ and $\constext{tur}$
(which by definition are always non-negative) correspond,
respectively, to the laminar and turbulent friction factors
\citep{
  WylieStreeter:78:book:Fluid,
  StreeterWylieBedford:98:book:Fluid,
  GerbeauPerthame:01:article:Derivation,
  LevermoreSammartino:01:article:A-shallow,
  Marche:07:article:Derivation
}.  The ground may also, due to
porosity, absorb water (by \emph{infiltration}) or inject water
(through \emph{recharge}) from and into the bulk. This mechanism is
modelled with the following \emph{permeable boundary condition}:
\begin{equation}
  \label{eq:permeable-Dirichlet-boundary-condition}
  \uu(t,x,z)
  \inner
  \normal[\W](x,z)
  =
  I(t, x)
  \quad \text{on } \mathcal B,
\end{equation}
where the infiltration function $I$ models the amount of water that
leaves ($I > 0$) or enters ($I < 0$) the flow per elementary boundary
element.

The term $\fricI$ models the friction effect that occurs when water that
is recharging through the ground (at average microscopic velocity rate
zero) connects with the flow. The magnitude of this effect is given by
the parameter $\alpha$, and hence our \indexemph{infiltration
mixing friction} law is given by
\begin{equation}
  \label{eq:def:infiltration-friction-law}
  \fricI
  :=
  \alpha\negpart I
  =
  \alpha\max(0,-I)
  .
\end{equation}
We note that the recharge-induced friction only occurs when water is
entering the flow (i.e. $I<0$), and is zero otherwise. Although $I$
should in principle be thought of as a function of $h$, $\uu$, and
possibly their derivatives---particularly $\mat\sigma[\uu]$, as in the
recognised Beavers--Joseph--Saffman model described, for instance, by
\citet{
  BeaversJoseph:67:article:Boundary,
  Saffman:71:article:On-the-Boundary,
  JagerMikelic:00:article:On-the-interface,
  BadeaDiscacciatiQuarteroni:10:article:Numerical
}---we ignore this in this paper and for simplicity consider the
function $I$ to be a given piecewise linear function of space-time.


We define $\W$'s tangential and outward unit normal vectors on $\mathcal B$ by
\begin{equation}
  \tangentialto\W(x,Z(x)) 
  = 
  \frac{\colvectwo{- 1}{- \pdx Z(x)}}{\sqrt{1 + \abs{\pdx Z(x)}^2}}
\end{equation}
and
\begin{equation}
  \normal[\W](x,Z(x)) 
  = 
  \frac{\colvectwo{\pdx Z(x)}{-1}}{\sqrt{1 + \abs{\pdx Z(x)}^2}},
\end{equation}
respectively, following the convention that the outward normal is the
tangential vector rotated by $\pi/2$ counterclockwise. It thus follows
that \eqref{eq:Navier-boundary-condition} and
\eqref{eq:permeable-Dirichlet-boundary-condition} on $\mathcal B$ can
be rewritten, respectively, as
\begin{gather}
  \label{eq:Navier-boundary-condition-coordinatewise}
  \begin{split}
    \frac{\mu\qp{\pd xv+\pd zu}
    \qpbig{1-\abs{\pd xZ}^2} 
    -
    2\mu\qp{\pd xu-\pd zv}{\pd xZ}}{\qpBig{1+\powabs2{\pd xZ}}^{\fracl12}}
    \\
    =
    \rho_0 \qp{k(u,v) + \fricI}\qp{u+v\pd xZ}
  \end{split}
  \intertext{and}
  \label{eq:permeable-Dirichlet-boundary-condition-coordinatewise}
  \vv
  -
  u
  \pdx Z(x)
  + I
  \sqrt{1 + \powabs2{\pdx Z
    }
  }
  =
  0
  .
\end{gather}

\subsection{The free surface}
\label{sec:free-boundary}
On the free surface, we neglect all other meteorological phenomena (such as evaporation) and consider only the addition of water in the form of direct rainfall and runoff. Assuming a \emph{kinematic boundary condition}, we set
\begin{equation}
\label{eq:kinematic-boundary-condition}
  \uu \inner \normal[\W] 
  = 
  \frac{\pdt H - R}{\sqrt{1+ \abs{\pdx H}^2}} 
  \quad \text{on } \mathcal F,
\end{equation}
where $R(t, x)$ is the recharge rate due to rainfall.  The unit tangential and normal vectors $\tangentialto \W$ and $\normal[\W]$ to the free surface can be explicitly computed in terms of $H$ as
\begin{equation}
  \tangentialto \W(x,H(t,x)) 
  = 
  \frac{\colvectwo{1}{\pdx H(t,x)}}{\sqrt{1 + \abs{\pdx H(t,x)}^2}}
\end{equation}
and
\begin{equation}
  \normal[\W](x,H(t,x)) 
  = 
  \frac{\colvectwo{-\pdx H(t,x)}1}{\sqrt{1 + \abs{\pdx H(t,x)}^2}},
\end{equation}
which leads to the following explicit form of \eqref{eq:kinematic-boundary-condition}:
\begin{equation}
  \label{eqn:kinematic-bc-on-top}
  \pdt H + u \pdx H - v = R \text{ on } \mathcal F.
\end{equation}
We also assume a stress condition on the free surface, given by
\begin{equation}
\qp{\mat\sigma\qb{\uu}\normal[\W]} 
  \inner 
  \tangentialto \W 
  = 
  -\rho_0 \fricR\uu
  \inner 
  \tangentialto \W,
\end{equation}
where we use the \indexemph{surface mixing friction} law
\begin{equation}
  \label{eq:mixing-friction-R}
  \fricR = \alpha R
\end{equation}
which takes into account the frictional effect of the additional
entering water's mixing with the flow from various sources (runoff,
direct rainfall, small-scale tributaries, for example) with $\alpha$
again representing the magnitude of this effect (see
\S\ref{subsec:friction-alpha} for more on mixing friction and
references). Using the tangential and normal vectors as above, this
condition becomes
\begin{equation}
  \label{bcfs1}
  \frac{\mu\qp{\pd xv+\pd zu}
    \qp{1-\powabs2{\pd xH}} 
    -
    2\mu\qp{\pd xu-\pd zv}{\pd xH}}{\sqrt{1+\powabs2{\pd xH}}}
  =
  -\rho_0 \fricR\qp{u+v\pd xH}
  .
\end{equation}
\section{Saint-Venant system with recharge via vertical averaging}
\label{SectionAveragedModel}
We now proceed to write the Navier--Stokes equations with adapted
boundary conditions in non-dimensional form. Under an assumption on
the shallowness of the ratio of the water height to the horizontal
domain (represented by a small parameter $\varepsilon$), we
formally make an asymptotic expansion of the Navier--Stokes system to
the hydrostatic approximation at first order. Finally, we derive the
Saint-Venant system through an integration on the water height.

This approach follows one established by
\citet{GerbeauPerthame:01:article:Derivation}, also found in
\citet{Ersoy:15:inproceedings:Dimension}, which we differ from in the
boundary conditions for Navier--Stokes that turn into different source
terms in the Saint-Venant's equation. Furthermore, we have to take
extra care in how we non-dimensionalise our additional precipitation,
infiltration, and friction terms. For simplicity, we will start from
the two-dimensional Navier--Stokes equations and obtain the
one-dimensional Saint-Venant's equations, although this procedure can
be employed to derive a two-dimensional analogue from the three
dimensional Navier--Stokes provided the boundary conditions are
modified accordingly.
\subsection{Dimensionless Navier--Stokes equations}
\label{sec:dimensionless-Navier-Stokes}
To derive the Saint-Venant model, we assume that the water height is small with respect to the horizontal length of the domain and that vertical variations in velocity are small compared to the horizontal variations. This is achieved by postulating a \indexemph{small parameter} ratio
\begin{equation}\index{$\varepsilon$}
  \varepsilon:= \frac{D}{L} = \frac{V}{U} \ll 1,
\end{equation}
where $D, L, V$, and $U$ are the \indexemph{scales} (or
\indexemph{units}) of, respectively, water height, domain length,
vertical fluid velocity, and horizontal fluid velocity. As a
consequence the time scale $T$ is such that
\begin{equation}
  T = \frac{L}{U}=\frac DV.
\end{equation}
We also choose the pressure scale to be
\begin{equation}
  \label{eq:def-pressure-scale}
  P:=\rho_0 U^2.
\end{equation}
The rationale for this choice is that we are focusing on the effect of the horizontal forces as mass per horizontal acceleration which has a force scale of
\begin{equation}
  F:=(DL^{2-1}\rho_0)(U\inverse T),
\end{equation}
and these forces are applied to vertical boundary scale to give the pressure scale
\begin{equation}
  F\inverseof{DL^{2-2}}
  =
  DL\rho_0U\inverse T\inverse D
  =
  \rho_0UL\inverse T
  =
  \rho_0U^2.
\end{equation}
It is convenient to define the spatial characteristic length, $L$, and horizontal velocity, $U$, (and by definition $T$) as finite constants with respect to $\ep\to 0$, while the water height and vertical velocity are defined as $D=\ep L$ and $V=\ep U$, respectively. This allows us to introduce the dimensionless quantities of time $\tilde t$, space $(\tilde{x}, \tilde{z})$, pressure $\tilde p$, and velocity field $(\tilde{u},\tilde{v})$ via the following scaling relations:
\begin{equation}
\label{eq:scaling-relations}
\begin{rcases}
\begin{dcases}
\begin{aligned}
\tilde t  &:= \frac{t}{T}, \quad & \tilde p(\tilde x, \tilde t, \tilde z) &:= \frac {p(x, t, z)}P,\\
\tilde x &:= \frac xL, \quad & \tilde u(\tilde x, \tilde t, \tilde z) &:= \frac {u(x, t, z)}U,\\
\tilde z &:= \frac{z}{D} = \frac z{\varep L}, \quad & \tilde v(\tilde x, \tilde t, \tilde z) &:= \frac{v(x, t, z)}{V} = \frac {v(x, t, z)}{\varep U} 
\end{aligned}
\end{dcases}
\end{rcases}.
\end{equation}
We also rescale the laminar and turbulent friction factors as, respectively,
\begin{equation}
  \label{eq:friction-factors-rescaled}
  \constext[0]{lam} 
  := \frac{\constext{lam}}{V} 
  = \frac{\constext{lam}}{\ep U}, 
  \quad 
  \constext[0]{tur} 
  := 
  \frac{\constext{tur}}{\ep},
\end{equation}
and the infiltration and rainfall rates as, respectively,
\begin{equation}
\label{eq:infiltration-rainfall-rescaled}
  \tilde I(\tilde t,\tilde x)
  := 
  \frac{I(t,x)}{V}, 
  \quad 
  \tilde R(\tilde t,\tilde x)
  := 
  \frac{R(t,x)}{V}
  .
\end{equation}
Note that in the assumed asymptotic setting, $\constext[0]{lam}$ and $\constext[0]{tur}$ are constants with respect to $\ep$, thus implying that $\constext{lam}$ and $\constext{tur}$ vanish linearly with $\ep\to0$. Finally, we define the following non-dimensional numbers:
\begin{equation}
\begin{aligned}
&\text{Froude's number,} &\Frhor &:= {U}/{\sqrt{\gravity D}},\\
&\text{Reynolds's number with respect to $\mu$,} &\Reynolds &:= {\rho_0 U L}/{\mu},
\end{aligned}
\end{equation}
and consider the following asymptotic setting
\begin{equation}
  \label{eq:Reynolds-asymptotic-setting}
  \Reynolds^{-1} = \varepsilon \mu_0,
\end{equation}
where $\mu_0$ is the \indexemph{viscosity}.

Using these dimensionless variables in the Navier--Stokes equations
\eqref{eq:incompressible-2d-Navier--Stokes} and
\eqref{NSTotalStressTensor},
and reordering the terms with respect to powers of $\varepsilon$, the
dimensionless incompressible Navier--Stokes system reads as
follows:\footnote{We bind all the ``tilde'' variables together, i.e.,
  $\tilde u$ is a function of $\tilde t,\tilde x,\tilde z$.  Hence
  variableless operators change accordingly, e.g., $\div\tildevec u$
  means $\div_{(\tilde x,\tilde z)}(\tilde u,\tilde v)$ when $\div\vec
  u$ means $\div_{(x,z)}(u,v)$.}
\begin{equation}
\begin{aligned}
  \label{AdimNSRearranged0}
  \divof{
  \tildevec u
  }
  &= 0 
  \\
  \pd{\tilde t} {\tilde u}
  + 
  \pd{\tilde x}\qb{\tilde u^2} 
  + 
  \pd{\tilde z}\qb{\tilde u \tilde\vv} 
  +  
  \pd{\tilde x}{\tilde p} 
  &=  
  \pd{\tilde z}\qb{\frac{\mu_0}{\varepsilon}\pd{\tilde z}\tilde u}
  + 
  \funkref[\varepsilon,\tildevec u]{remainder1}{\varrho}
  \\
  \pd{\tilde z}\tilde p 
  &= 
  -\frac{1}{\Frhor^2} 
  + 
  \funkref[\varepsilon,\tildevec u]{remainder2}{\varrho}
\end{aligned}
\end{equation}
where 
\begin{equation}
  \funkdef[\varepsilon,\tildevec u]{remainder1}{\varrho}
  := 
    \varepsilon\mu_0\Big(2
    \pd{\tilde u}{\tilde x\tilde x}
    + \pd{\tilde v}{\tilde z\tilde x}\Big)
\end{equation}
and
\begin{equation}
  \funkdef[\varepsilon,\tildevec u]{remainder2}{\varrho}
  :=  
    \varepsilon\mu_0
    \Big(\pd{\tilde x\tilde z}{\tilde u}
    +
    \varepsilon^2
    {\pd{{\tilde x}\tilde x}{\tilde v}}
    +
    2\pd{\tilde z\tilde z}{\tilde \vv}\Big)
  - 
  \varepsilon^2
  \Big(
    \pd{\tilde t}{\tilde\vv}
    +\pd{\tilde x}\qb{{\tilde u\tilde \vv}}
    +\pd{\tilde z}\qb{\tilde \vv^2}
  \Big).
\end{equation}
Assuming $\tildevec u$ has bounded second derivatives, definitions
\eqref{remainder1} and \eqref{remainder2} formally lead to
\begin{equation}
  \funkref[\varepsilon,\tildevec u]{remainder1}{\varrho} = \Oh(\varepsilon) \quad \text{and} \quad
  \funkref[\varepsilon,\tildevec u]{remainder2}{\varrho} = \Oh(\varepsilon).
\end{equation}

On the wet boundary $\mathcal B$, recalling the scaling relations
\eqref{eq:scaling-relations} and
\eqref{eq:infiltration-rainfall-rescaled}, and noting that
\begin{equation}
  \begin{aligned}
    \frac{\partial Z}{\partial x}
    =
    \frac{\varep L}{L} \frac{\partial \tilde Z}{\partial \tilde x}
    =
    \varep \partial_{\tilde x} \tilde Z
    ,
  \end{aligned}
\end{equation}
the dimensionless Navier boundary 
condition \eqref{eq:Navier-boundary-condition-coordinatewise} implies
\begin{equation}
  \begin{aligned}
    \qb{
      \frac{\partial_{\tilde{z}}\tilde{u}}{\varep \Reynolds}}_{\mathcal{B}}
    = &
    \qp{
      \frac{C_{\text{lam}}}{U}\tilde{u}
      +
      C_{\text{tur}}\qp{\abs{\tilde u}
        + \varep\abs{\tilde v}
      }\tilde{u}
      +
      \varep \fricI[\tilde{I}]\tilde u
    }
    \frac{
      \sqrt{1 + \varep^2(\partial_{\tilde{x}} \tilde{Z})^2}
    }{
      1 - \varep^2(\partial_{\tilde{x}} \tilde{Z})^2}
    \\
    &+\underbrace{
      \varep^2\partial_{\tilde{x}} \tilde{Z}
      \qp{
        \frac{C_{\text{lam}}}{U}\tilde{v}
        +
        C_{\text{tur}}\qp{\abs{\tilde u}
          + \varep\abs{\tilde v}}\tilde{v}
        + \varep \fricI[\tilde{I}]\tilde v}
      \frac{
        \sqrt{1 + \varep^2(\partial_{\tilde{x}} \tilde{Z})^2}
      }{
        1 - \varep^2(\partial_{\tilde{x}} \tilde{Z})^2}
    }_{\BigO{\varep^2}}
    \\
    &-
    \underbrace{
      \frac{\varep}{\Reynolds}
      \qp{
        \partial_{\tilde{x}}\tilde{v}
        + \frac{2\partial_{\tilde{x}} \tilde{Z}
          \qp{
            \partial_{\tilde{z}}\tilde{v}
            - \partial_{\tilde{x}}\tilde{u}
          }
        }{
          1 - \varep^2(\partial_{\tilde{x}} \tilde{Z})^2}
      }
    }_{\BigO{\sfrac{\varep}{\Reynolds}}}.
  \end{aligned}
\end{equation}
Applying the non-dimensional friction factors
\eqref{eq:friction-factors-rescaled} and recalling
\eqref{eq:Reynolds-asymptotic-setting}, we get
\begin{equation}
  \begin{aligned}
    \label{eq:Navier-boundary-condition-adim-rearranged}
    \qb{\frac{\partial_{\tilde{z}}\tilde{u}}{\varep \Reynolds}}_{\mathcal{B}}
    &=
    \varep\qp{C_{\text{lam},\varep}\tilde{u} + C_{\text{tur},\varep}\qp{\abs{\tilde u}
        + \varep\abs{\tilde v}}\tilde{u}
      + \fricI[\tilde{I}]\tilde u}\frac{\sqrt{1
        + \varep^2(\partial_{\tilde{x}} \tilde{Z})^2}}{1 - \varep^2(\partial_{\tilde{x}} \tilde{Z})^2} + \BigO{\varep^2} \\
    &=
    \varep
  \qp{
    C_{\text{lam},\varep}\tilde{u}
    + C_{\text{tur},\varep}\abs{\tilde{u}}\tilde{u}
    + \fricI[\tilde{I}]\tilde u
  }
  + \BigO{\varep^2}
  \\
  &=
  \varep \qp{k_0(\tilde{u}) + \fricI[\tilde{I}]}\tilde{u} + \BigO{\varep^2},
\end{aligned}
\end{equation}
with \indexemph{asymptotic friction laws}
\begin{equation}
  \begin{aligned}
    \label{eqn:wet-boundary-friction-nondim}
    k_0(\xi) &:= C_{\text{lam},0} + C_{\text{tur},0} \abs{\xi} \text{ for } \xi \in \reals\\
    \fricI[\tilde I]&:= \alpha\negpart I=-\min(0, \tilde{I})
    \text{ for } \alpha \in \reals,
  \end{aligned}
\end{equation}
on the wet boundary. The permeable boundary condition \eqref{eq:permeable-Dirichlet-boundary-condition-coordinatewise} reads
\begin{equation}
  \label{eq:permeable-Dirichlet-boundary-condition-adim-rearranged}
   \tilde v 
   =
   \tilde  u \pdx  Z 
   - I \sqrt{1+\varepsilon^2(\pdx  Z)^2} 
   = 
   \tilde u \pdx  Z 
   - I + \Oh(\varepsilon) .
\end{equation}
For the boundary conditions on the free surface $\mathcal{F}$, applying our non-dimensionalisation approach to the kinematic boundary condition \eqref{eqn:kinematic-bc-on-top} we derive
\begin{equation}
\label{eqn:kinematic-bc-nondimensional}
\partial_{\tilde t} \tilde H + u \partial_{\tilde x} \tilde H - \tilde v = \tilde R,
\end{equation}
while we can non-dimensionalise \eqref{bcfs1} in the same manner as the Navier boundary condition \eqref{eq:Navier-boundary-condition-coordinatewise} on the wet boundary $\mathcal{B}$, giving
\begin{equation}
\begin{aligned}
  \label{eq:free-surface-nondimensional}
  \qb{\frac{\partial_{\tilde{z}}\tilde{u}}{\varep \Reynolds}}_{\mathcal{F}}
  &
  =
  -\varep \fricR[\tilde{R}]\tilde{u} + \BigO{\varep^2},
\end{aligned}
\end{equation}
with free surface asymptotic friction law
\begin{equation}
  \begin{aligned}
    \label{eqn:free-surface-friction-nondim}
    \fricR[\tilde R]
    =
    \alpha \tilde R \text{ for } \alpha \in \mathbb{R}.
  \end{aligned}
\end{equation}

\subsection{First order approximation of the dimensionless Navier--Stokes equations}
Dropping all terms of $\Oh(\varepsilon)$ and above in equations \eqref{AdimNSRearranged0}-\eqref{eq:free-surface-nondimensional}, we deduce the hydrostatic approximation of the dimensionless Navier--Stokes system
\begin{align}
  \pdx  {\ue} + \pd z\ve &= 0 
  \label{HS0}
  \\
  \pdt  {\ue} + \pd{x}\qb{{\ue}^2}
  + \pd[{\ue} \ve]z 
  + \pdx  \pe 
  &
  = 
  \pd z 
  \left[{\frac{\mu_0}{\varepsilon}\pd z {\ue}}\right] 
  \label{HS1}
  \\
  \pd z\pe &= - \frac{1}{\Frhor^2} \label{HS2},
\end{align}
whilst the boundary conditions, as a result of the asymptotic setting \eqref{eq:Reynolds-asymptotic-setting}, simplify to 
\begin{equation}
\begin{aligned}
\label{eqn:wet-boundary-nondim-first-order}
\qb{\frac{\mu_0}{\varepsilon}\pdz\ue} = \qp{k_0(\ue) + \fricI}\ue \quad \text{and} \quad \ve = \ue \pdx Z - I \text{ on } \mathcal{B},
\end{aligned}
\end{equation}
and
\begin{equation}
\begin{aligned}
\label{eqn:free-surface-nondim-first-order}
\qb{\frac{\mu_0}{\varepsilon}\pdz\ue} = - \fricR\ue \quad \text{and} \quad \pdt H + \ue \pdx H - \ve = R \text{ on } \mathcal{F},
\end{aligned}
\end{equation}
in view of equations
\eqref{eq:Navier-boundary-condition-adim-rearranged},
\eqref{eq:permeable-Dirichlet-boundary-condition-adim-rearranged}, \eqref{eq:free-surface-nondimensional}, and \eqref{eqn:kinematic-bc-nondimensional}, respectively. Here, $(\ue, \ve, \pe)$ represents the solution of the first-order dimensionless Navier-Stokes system.

Vertically integrating both members of equation \eqref{HS2} over $[z,H(t,x)]$, we obtain the hydrostatic pressure
\begin{equation}
  \label{p}
  \pe(t,x,H) - \pe(t,x,z) 
  =   
  - \frac{1}{\Frhor^2}(H(t, x)-z)
  .
\end{equation}
Assuming that the pressure exerted by the rain on the free surface $\pe(t, x, H) = p_c$ for some constant $p_c \in \reals$ (as we neglected all other meteorological phenomena), this becomes
\begin{equation}
  \label{p}
  \pe(t,x,z) 
  = \frac{1}{\Frhor^2}(H(t, x)-z) + p_c
  .
\end{equation}

Moreover, identifying terms at order $\sfrac{1}{\varepsilon}$ in
\eqref{HS1}, \eqref{eqn:wet-boundary-nondim-first-order}, and \eqref{eqn:free-surface-nondim-first-order}, we obtain the
\indexemph{motion by slices} decomposition
\begin{equation}
  \label{MotionBySlices}
    {\ue}(t,x,z)   
    = 
    u_0(t, x) 
    +  \Oh(\varepsilon)  
\end{equation}
for some function $u_0 =u_0(t, x)$, as a consequence of
\begin{gather}
  \pdz\qb{\mu_0 \pd z {\ue} } =  \Oh(\varepsilon), \text{ for } z \in (Z(x),H(t, x)) 
  \intertext{with} 
  \qb{\mu_0 \pd z {\ue} }_{| z = Z(x)} =  \Oh(\varepsilon)  
  \tand   
  \qb{\mu_0 \pd z {\ue} }_{| z = H(t, x)} =  \Oh(\varepsilon) .
\end{gather}

Noting $\average{{\ue(t,x)}}$ as the mean speed of the fluid over the section $[Z(x),H(t, x)]$,
\begin{equation}
  \average{{\ue}(t,x)} 
  = 
  \frac{1}{h(t, x)}\int_{Z(x)}^{H(t, x)} {\ue}(t,x,z) \d z,
\end{equation}
we are able to use the following approximations and drop the
first and higher order terms in $\ep$:
\begin{equation}
  \label{FirstOrder1}
  \ue(t, x, z) 
  =
  \average{{\ue}(t, x)}
  +
  \Oh(\varepsilon)
  \tand
  \average{{\ue(t,x)}^2}
  =
  \average{\ue(t,x)}^2 
  +
  \Oh(\varepsilon)
  .
\end{equation}

\subsection{The Saint-Venant system with recharge}
Keeping in mind \eqref{FirstOrder1} and integrating the indicator transport equation \eqref{eqn:indicator-transport} for $z \in [Z(x),H(t, x)]$, we get
\begin{equation}
  \begin{split}
    0 
    &= 
    \int_{Z(x)}^{H(t, x)} 
    \pdt  \Phi(t,x,z) 
    + 
    \pdx \qb{\Phi{\ue}}
    + 
    \pd z\qb{
      \Phi
      \ve}
    \d z 
    \\
    &= 
    \pdt  h + \pdx  q - \qb{\pdt  H + \ue \pdx  H - \ve}_{z=H(t, x)} ]
      + \qb{\ue \pdx  Z -\ve}_{z=Z(x)},
  \end{split}
\end{equation}
where $q$ is the discharge defined by
\begin{equation}
  q(t,x):=
  \average{\ue(t,x)}h(t,x).
\end{equation} 
In view of the penetration condition
\eqref{eqn:wet-boundary-nondim-first-order} and the kinematic boundary condition \eqref{eqn:free-surface-nondim-first-order}, we
deduce the \indexemph{mass-balance} equation:
\begin{equation}
  \label{SLMasse}
  \pdt  h + \pdx  q = S
\end{equation}
where the source term $S := R - I$ measures the gain or loss of water
through the (nonnegative) \indexemph{recharge rate} $R$ (ultimately
from rainfall) and (signed) \indexemph{infiltration rate},
respectively.

Keeping equations \eqref{p}, \eqref{MotionBySlices}, and
\eqref{FirstOrder1} in mind and thanks to the penetration condition
\eqref{eqn:wet-boundary-nondim-first-order} and the kinematic boundary condition \eqref{eqn:free-surface-nondim-first-order}, integrating the left-hand side of \eqref{HS1} for $z \in [Z(x),H(t, x)]$, we get
\begin{equation}
  \label{SLMomLHS}
  \begin{aligned}
    \int_{Z(x)}^{H(t,x)}
    \LHS{HS1}
    \d z
    &= 
    \pdt  q + \pdx \qb{\frac{q^2}{h}+ \frac{h^2}{2 \Frhor^2}} 
    + 
    \frac{h}{\Frhor^2}\pdx  Z
    \\
    &
    \phantom=
    - 
    \evalat{\qp{\pdt  H + {\ue} \pdx  H - \ve} {\ue}}{(t,x,H(t, x))}
    \\
    &
    \phantom=
    +
    \evalat{\qp{{\ue} \pdx  Z - \ve} {\ue}}{(t,x,Z(x))}
    \\
    &= 
    \pdt  q 
    + 
    \pdx \qb{\frac{q^2}{h}+ \frac{h^2}{2 \Frhor^2} } 
    + 
    \frac{h}{\Frhor^2}\pdx  Z 
    \\
    &
    \phantom=
    - 
    R \evalat{\ue}{(t,x,H(t, x))}
    +
    I \evalat{\ue}{(t,x,Z(x))}
    \\
    &= 
    \pdt  q 
    + 
    \pdx \qb{\frac{q^2}{h}+ \frac{h^2}{2 \Frhor^2} }
    + 
    \frac{h}{\Frhor^2}\pdx  Z - S \frac{q}{h},
  \end{aligned}
\end{equation}
where $S$ is again defined as above. Now, integrating the right-hand side of \eqref{HS1} for $z
\in [Z(x),H(t, x)]$ using the wet boundary condition \eqref{eqn:wet-boundary-nondim-first-order} and
the free surface boundary condition \eqref{eqn:free-surface-nondim-first-order}, we obtain:
\begin{equation}\label{SLMomRHS}
\begin{aligned}
    \int_{Z(x)}^{H(t,x)}\RHS{HS1}\d z 
    &= 
    \evalat{\frac{\mu_0}{\varepsilon} \pd z {\ue} }{z=H(t, x)}
    -
    \evalat{\frac{\mu_0}{\varepsilon} \pd z {\ue} }{z=Z(x)}\\
    &= 
    - \qp{\fricR + \fricI + k_0\qp{\frac{q}{h}}} \frac{q}{h},
\end{aligned}
\end{equation}
where the friction factors $\fricR$, $\fricI$, and $k_0$ are defined
by formulas \eqref{eqn:free-surface-friction-nondim} and
\eqref{eqn:wet-boundary-friction-nondim}, respectively. Finally,
multiplying both sides of each of \eqref{SLMasse}, \eqref{SLMomLHS},
and \eqref{SLMomRHS} by $\rho_0 U^2 /D$, and recalling the
mass-balance (\ref{SLMasse}), we obtain the following
\indexemph{Saint-Venant system with recharge}:
\begin{equation}
  \begin{aligned}
    \label{ModelSL}
    \pdt  h +\pdx  q
    &= S:=R-I,
    \\
    \pdt  q
    +\pdx \qb{\frac{q^2}{h}
      +\gravity \frac{h^2}{2}
    }
    &=
    - \gravity   h \pdx  Z 
    + S \frac{q}{h}
    - \qp{\fricR + \fricI + k_0\qp{\frac{q}{h}}}
    \frac{q}{h}
    \\
    \text{where }
    q&=hu
    ,
  \end{aligned}
\end{equation}
which we study in the rest of this paper.
\begin{Exa}[lake at rest and filling the lake]
  \label{item:the:main:steady-state}
  The still water steady state \akaindexemph{lake at rest} reads
  \begin{equation}
    \label{SteadyState}
    q\ideq u\ideq S\ideq 0 \hbox{ and } h + Z\ideq H_0 \hbox{ for some constant  } H_0 > 0.
  \end{equation}
  This is a classical example used in testing the conservation properties of numerical
  schemes. Note that an interesting nontrivial example for numerical tests
  is $S=R-I\ideq0$ with $R=I>0$.

  Another simple (yet important) example is a space-time uniform
  filling of a lake with initial height $h(0,\cdot)\ideq H_0$ as
  in~(\ref{SteadyState}) with a constant time-space $S>0$ and a
  spatially constant momentum with periodic boundary conditions.  In
  this case, symmetry implies that system (\ref{ModelSL}) simplifies
  to
  \begin{equation}
    \pdt h\ideq S
    \and
    \pdt q\ideq
    S\frac qh
    -\qp{\fricR+\fricI+k_0\qp{\frac qh}}\frac qh,
  \end{equation}
  with given initial conditions.
\end{Exa}
\subsection{Mixing friction}
\label{subsec:friction-alpha}
In adapting the boundary conditions of the Navier-Stokes equations in
\S\ref{sec:wet-boundary} and \S\ref{sec:free-boundary}, we included
additional terms $\fricR$ and $\fricI$ that model the friction effect
that occurs when water that is falling on the free surface or
recharging through the ground, respectively, connects with the
flow. The inclusion of these terms avoids certain paradoxical
outcomes, such as perpetual motion, that otherwise occur when the
terms are omitted.  Such terms arise naturally from microscopic
effects and have been discussed in the hydrology literature, for
instance, by \citet{
  Wenzel:70:techreport:The-Effect
  ,YoonWenzel:71:article:Mechanics
  ,ShenLi:73:article:Rainfall
  ,LuChenChangLu:98:article:Characteristics
} where the laws are empirically derived from measurements.  While
the friction can be quite complex in behaviour depending on many of of
the involved quantities, but they share the characteristics of being
monotone, zero at zero, and possibly homogeneous in the velocity $u$
and the rates of additional water $R$ and $I$, we consider the simplest
such behaviour by taking.

To see the influence of the friction effect $\alpha$ on the solution,
we consider an idealised scenario that will enable us to calculate an
\emph{exact solution} to system \eqref{ModelSL}. We take a constant
rainfall--runoff process on a river spanning spatial domain
$x\in[0,10]$ and time domain $t\in[0,1]$, with topography $Z\ideq0.1$;
for simplicity, we assume the infiltration $I\ideq0$ and a linear
recharge friction $\fricR=\alpha{R}$. We prescribe periodic boundary
conditions and assume a constant initial height and discharge of
\begin{equation}
  \begin{aligned}
    h(0, x) = q(0, x) = 1\text{ for }0<x<10.
  \end{aligned}
\end{equation}
The rainfall intensity is applied uniformly on the river as a function
of time up to the final time $T=1$:
\begin{equation}
  \begin{aligned}
    R(t) = 1.0\text{ for }0<t<1.
  \end{aligned}
\end{equation}
\begin{figure}[ht]
  \centering
  \scalebox{.8}{\begin{tikzpicture}
      \begin{scope}
        \clip (0,1.2) -- (10,1.2) -- (10,4.5) -- (0,4.5);
        \fill [cyan,opacity=0.3] (-0.5,-0.5) rectangle (10,5);
      \end{scope}
      
      \begin{scope}
        \clip (0,0) -- (0,1.2) -- (10,1.2) -- (10,0);
        \fill [brown, opacity=0.3] (-0.5,-0.5) rectangle (10,4);
      \end{scope}
      
      \draw [cyan, thick] plot[smooth, domain=0:10]({\x}, {4.5});
      \draw [brown, thick] plot[smooth, domain=0:10]({\x}, {1.2});
      \draw [black, -latex] (-0.5,0) -- (10.5,0);
      \draw [black, -latex] (0,0) -- (0,5.5);
      
      \draw (0,0) node [below] {0};
      \draw (10,0) node [below] {10};
      
      \draw [cyan, ultra thick, -latex] (4,2.85) -- node [above] {$q(x,0) = 1$} (8,2.85);
      \draw [cyan, ultra thick, latex-latex] (2,1.2) -- node [above, rotate=90] {$h(x,0) = 1$} (2,4.5);
      \draw [black, thick, -latex] (5,5.5) node [above] {$R(t) = 1$} -- (5,4.5);
      
      \draw [cyan, ultra thick, domain=90:270, latex-] plot ({10.25-cos(\x)}, {2.85+0.5*sin(\x)});
      \draw [cyan, ultra thick, domain=90:270, -latex] plot ({-0.25+cos(\x)}, {2.85+0.5*sin(\x)});
  \end{tikzpicture}}
  \caption{Diagram of the idealised scenario we are considering.}
\end{figure}
\par
Since the rain function, initial height, and initial discharge do not have any dependence on $x$, we can ignore the spatial derivatives in both the mass and momentum equations. Under these assumptions, \eqref{ModelSL} simplifies to:
\begin{equation}
  \begin{aligned}
    \pdt h &= 1,\\
    \pdt q &= \qp{1 - \alpha} \frac{q}{h}\quad \text{where } q = hu.
  \end{aligned}
\end{equation}
We can solve explicitly for $h$ and $q$ (which are constant in space)
and the corresponding velocity
\begin{equation}
    h(t,x) 
    = t + 1,
    \quad
    q(t,x) 
    = (t + 1)^{1 - \alpha},
    \quad
    u(t,x) 
    =
    (t + 1)^{-\alpha}.
\end{equation}
Using these equations, we can plot how the \emph{momentum} and
\emph{velocity} change over time depending on the mixing friction
coefficient $\alpha$. We consider three cases, reported in
Figure~\ref{fig:momentum-paradox-exact-solution-plots},
\begin{itemize}
\item[${\alpha< 1}$]
  water entering the course produces
  too little friction to slow the flow down.  In this regime, the
  velocity remains constant or decreases over time while the \emph{momentum
    increases} at the rate of increase in water height. The exclusion, or
  weakness, of the additional friction terms $\fricI$ and $\fricR$
  leads to a paradox, as momentum is being produced even though the
  rain is assumed to fall with zero momentum. We can conclude,
  therefore, that this regime leads to an unrealistic physical system.
\item[${\alpha=1}$]
  the friction generated by the rain
  causes a reduction in velocity that balances the increase in
  height. We are thus in the situation where the additional terms in
  the conservation of momentum equation cancel one another out.  In
  this regime, momentum is conserved over time due to the balance
  between the increase in momentum from the $Ru$ term and the
  reduction in momentum from the additional friction term; this leads
  to a reduction in velocity over time. As the momentum is conserved
  within the system, this regime is physically realistic.
\item[${\alpha > 1}$]
  the friction generated by the rain
  exceeds the increase in momentum from the $Ru$ term. We are thus in
  the situation where the additional friction term is dominant.  In
  this regime, the additional friction term slows the flow at a faster
  rate than the increase in water height; therefore, momentum and
  velocity both decrease over time. The reduction in momentum due to
  the higher friction effect means the regime is physically realistic.
\end{itemize}

We can conclude from considering this idealised scenario that the
additional friction terms cannot be omitted, as doing so leads to a
paradoxical situation where the momentum increases over time even
though the rain is assumed to be added to the system with zero
momentum. We therefore only have physically realistic solutions when
the friction $\alpha \geq 1$, as it is within this regime that the
momentum is either conserved or decreasing. As is clear from the
equations, this occurs when the reduction in velocity is large enough
to either compensate or balance the increase in water height.
\begin{figure}
  \scalebox{0.6}{%
    \begin{tikzpicture}
      \begin{axis}[width=0.8\textwidth, height=0.55\textwidth, xlabel={$t$}, title=\text{momentum}, ylabel={$q(t,x)$}, legend pos=north west, legend cell align={left}, ymin=0, ymax=2.5, xmin=0, xmax=1]
        \addplot [a, very thick, domain=0:1, smooth, samples=100] {(x+1)^(1)};
        \addlegendentry{\ $\alpha = 0$}
        \addplot [b, very thick, domain=0:1, smooth, samples=100] {(x+1)^(0)};
        \addlegendentry{\ $\alpha = 1$}
        \addplot [c, very thick, domain=0:1, smooth, samples=100] {(x+1)^(-1)};
        \addlegendentry{\ $\alpha = 2$}
      \end{axis}
    \end{tikzpicture}
  }
  \scalebox{0.6}{%
    \begin{tikzpicture}
      \begin{axis}[width=0.8\textwidth, height=0.55\textwidth, xlabel={$t$}, title=\large{velocity}, ylabel={$u(t,x)$}, legend pos=north west, legend cell align={left}, ymin=0, ymax=2.5, xmin=0, xmax=1]
        \addplot [a, very thick, domain=0:1, smooth, samples=100] {(x+1)^(0)};
        \addlegendentry{\ $\alpha = 0$}
        \addplot [b, very thick, domain=0:1, smooth, samples=100] {(x+1)^(-1)};
        \addlegendentry{\ $\alpha = 1$}
        \addplot [c, very thick, domain=0:1, smooth, samples=100] {(x+1)^(-2)};
        \addlegendentry{\ $\alpha = 2$}
      \end{axis}
    \end{tikzpicture}
  }
  \caption{
    \label{fig:momentum-paradox-exact-solution-plots}
    Exact solution's momentum and velocity in three cases of mixing
    friction coefficient $\alpha$ as discussed in
    \S~\ref{subsec:friction-alpha}. The case $\alpha=0$ means no
    mixing friction whatsoever and leads to a \indexemph{physical
      paradox} where momentum increases in time as mass increases.
    The case $\alpha=1$ is critical and ensures no artificial momentum
    gain, but this case is also an unrealistic idealised situation
    because it implies momentum conservation in presence of mixing
    friction.  Finally the case $\alpha>1$ is the most physically
    relevant.}
\end{figure}
\section{Entropy}
\label{sec:entropy-condition}
\begin{LongChanges}
  Entropy plays an essential role in the analytical and numerical
  understanding of conservation and balance laws
  \citep{Evans:13:booklet:Entropy,Serre:99:book:Systems,Dafermos:10:book:Hyperbolic}.
  In this section we study the effect of the rainfall terms and
  associated mixing friction on the entropy--entropy-flux pair, for the
  Saint-Venant system \eqref{ModelSL}, in particular the effect of the
  additional rainfall terms have on entropy production.
\end{LongChanges}
\begin{The}[hyperbolicity and stability of~\eqref{ModelSL}]
  \label{THM_ESV}
  Let $(h,q)$ with satisfy the Saint-Venant system
  with recharge~\eqref{ModelSL} for a given topography $Z$,
  rainfall $R$ and infiltraion $I$, with velocity
  \begin{equation}
    u := \frac qh
  \end{equation}
  and \indexemph{total head} $\psi := \hat\psi(h,q,Z)$ defined by
  \begin{equation}
    \begin{aligned}
      \label{eqn:defn-entropy-head}
      \hat \psi(h, q, Z) := \frac{q^2}{2h^2} + \gravity\qp{ h + Z},
      \text{ for }(h,q,Z)\in\R+\times\R2.
    \end{aligned}
  \end{equation}
  Then the following hold:
  \begin{enumerate}[label=(\alph*)]
  \item
    System \eqref{ModelSL} is strictly hyperbolic on the set
    $\setof{h>0}$.\label{item:the:main:hyperbolicity}
  \item
    \label{item:the:main:velocity-head-relation} If $(h,q)$ is
    smooth and $h>0$, we have the \indexemph{velocity balance
      equation}
    \begin{equation}
      \label{TotalHead}
      \pdt u + \pdx \psi = - \frac{\qp{\fricR + \fricI +  k_0(u)}u}{h}.   
    \end{equation}
  \end{enumerate}
\end{The}
\begin{proof}
  The proof follows standard arguments from conservation laws
  and is provided for self-containment's sake.
  We consider each statement in turn.
  \begin{LetterList}
  \item
    The Jacobian of (\ref{ModelSL})'s \indexemph{flux function}
    $(q,h)\mapsto(q,\fracl{q^2}h+\gravity\fracl{h^2}2)$ is given by
    \begin{equation}
      \begin{aligned}
        \label{eq:Jacobian}
        \vect J=
        \hatmat{J}(h, q)
        :=
        \begin{bmatrix}
          0
          &
          1
          \\
          - \sfrac{q^2}{h^2} + \gravity h
          &
          \sfrac{2q}{h}
        \end{bmatrix},
      \end{aligned}
    \end{equation}
    with eigenvalues
    \begin{equation}
      \begin{aligned}
        \lambda_{1,2}
        =
        \hat\lambda_{1,2}(h, q) = \frac qh \pm \sqrt{\gravity h}
        .
      \end{aligned}
    \end{equation}
    For these eigenvalues to be real and distinct, we require
    that $h > 0$; the Jacobian matrix is thus diagonalizable and~
    system \eqref{ModelSL} is strictly hyperbolic on the set $\qc{h > 0}$.
  \item
    We rewrite the conservation of momentum equation in system
    \eqref{ModelSL} in terms of the unknowns $(h, u)$, with
    $u=\fracl{q}h$, as
    \begin{equation}
      \label{Momentum_Rewrite}
      \pdt \qb{hu} + \pdx \qb{hu^2 + \gravity \frac{h^2}{2}}
      =
      - \gravity h \pdx  Z
      + Su
      - \qp{\fricR + \fricI + k_0(u)}u.
    \end{equation}
    Applying the product rule to the first term of
    \eqref{Momentum_Rewrite} and substituting in the conservation of
    mass equation, we get
    \begin{multline}
      \label{Momentum_Reorder1}
      h \pdt u
      +
      u \qp{S - \pdx \qb{hu}}
      + \pdx \qb{hu^2} + \pdx \qb{\gravity\frac{h^2}{2}}
      \\
      =
      - \gravity h \pdx  Z + Su - \qp{\fricR + \fricI + k_0(u)}u,
    \end{multline}
    from which we can cancel $Su$ on both sides. Using the
    product rule again, we have that
    \begin{equation}
      u\pdx \qb{hu} = \pdx \qb{hu^2} - hu\pdx u,
    \end{equation}
    which can be substituted into \eqref{Momentum_Reorder1} to give
    \begin{multline}
      \label{Momentum_Reorder2}
      h \pdt u
      - \pdx \qb{hu^2}
      + hu \pdx u + \pdx \qb{hu^2}
      + \pdx \qb{\gravity\frac{h^2}{2}}
      \\
      =
      - \gravity h \pdx  Z - \qp{\fricR + \fricI + k_0(u)}u,
    \end{multline}
    enabling us to now cancel $\pdx \qb{hu^2}$. We note that
    \begin{equation}
      \pdx \qb{\gravity \frac{h^2}{2}} = h \pdx \qb{\gravity h}.
    \end{equation}
    Substituting this into \eqref{Momentum_Reorder2} and dividing by $h$
    throughout, we get
    \begin{equation}
      \pdt u
      + u \pdx u
      + \pdx[{\gravity h}]
      =
      -
      \gravity  \pdx  Z
      -
      \frac{\qp{\fricR + \fricI + k_0(u)}u}{h}.
    \end{equation}
    Making the further substitution
    \begin{equation}
      u \pdx u = \pdx[u^2]/2
    \end{equation}
    and grouping derivatives of $x$, we have
    \begin{equation}
      \pdt u + \pdx \psi
      =
      - \frac{\qp{\fricR + \fricI + k_0(u)}u}{h}
    \end{equation}
    as required.
  \end{LetterList}
  The theorem is thus proven.
\end{proof}
\begin{Obs}[friction effects]
  It is worth noting that the only way $S=R-I$ enters in the velocity
  balance equation ~(\ref{TotalHead}) is through the friction terms
  $\fricR$ and $\fricI$.  This is a further indication that these
  terms are necessary, especially for small water height $h$, which as
  a denominator of the right-hand side in (\ref{TotalHead}) amplifies
  the effects of friction as these occur in the Navier--Stokes on a
  layer close to the boundaries.
\end{Obs}
\begin{The}[entropy production]
  \label{THM_ESV2}
  Consider the \indexemph{entropy} \index{$\hat E$}
  and
  \indexemph{entropy flux} \index{$\hat\Psi$} respectively defined by
  \begin{gather}
    \label{eqn:defn-entropy}
    \hat E(h, q, Z)
  :=
    \frac{q^2}{2h}+\gravity h\qp{\frac h2+Z}
    =
    \frac{hu^2 + \gravity h^2}{2}  + \gravity hZ
    ,
    \\
    \label{eqn:def:entropy-flux}
    \hat \Psi(h, q, Z, S) :=
    \qp{\hat E(h,q) + \frac{\gravity h^2}{2}}\frac qh
    -\gravity\Theta
    .
  \end{gather}
  where
  \begin{equation}
    \Theta(t,x):=\int_0^xS(t,s)Z(s)\d s
    \Foreach
    t>0,x\in\reals
  \end{equation}
  The pair of functions $(\hat E, \hat \Psi)$
  forms a mathematical \indexen{entropy}--\indexen{entropy-flux} pair for system
  \eqref{ModelSL} when $S\ideq0$ and $k_i=0$ for $i=0,\pm$.
  Furthermore for smooth solutions $(h,q)$ of (\ref{ModelSL}),
  the functions 
  \begin{equation}
    E:=\hat E(h,q,Z)\tand\Psi:=\hat\Psi(h,q,Z,S)
  \end{equation}
  satisfy the following \indexemph{entropy production} relation
  \begin{equation}
    \begin{aligned}
      \label{eqn:entropy-condition-SWE}
      \partial_t E
      + \partial_x \Psi
      = S\qp{\frac{u^2}{2}
        + \gravity h
      }
      - \qp{\fricR + \fricI + k_0(u)}u^2.
    \end{aligned}
  \end{equation}
\end{The}
\begin{Obs}[entropy--entropy-flux pairs and entropy production]
  \label{obs:entropy-entropy-flux-and-production}
  The entropy--entropy-flux pair of Theorem~\ref{THM_ESV2}
  is typical for shallow water equations
  \citep[Th.~2.1, e.g.]{AudusseBristeauPerthame:00:techreport:Kinetic}.
  Indeed (\ref{eqn:entropy-condition-SWE}), when $S=0$ and $k_i=0$, for $i=0,\pm$, generalises
  the well-known (zero) \indexemph{entropy condition}
  \begin{equation}
    \partial_t E + \partial_x \Psi = 0,
  \end{equation}
  which implies that $({\hat E},{\hat\Psi})$ is an entropy--entropy-flux pair.
  The additional term $-\gravity\Theta$ in (\ref{eqn:def:entropy-flux}) corresponds
  to the flux of entropy due to the added or subtracted rain; thanks to this term
  the entropy production is frame invariant.
  
\end{Obs}
\begin{Proofof}{Theorem~\ref{THM_ESV2}}
  In view of \ref{obs:entropy-entropy-flux-and-production} it is
  enough to prove only (\ref{eqn:entropy-condition-SWE}).  We
  begin by recalling that we showed in
  Theorem~\ref{THM_ESV}\ref{item:the:main:velocity-head-relation} that
  the conservation of momentum equation can be rewritten in terms of
  the velocity $u$ as
  \begin{equation}
    \begin{aligned}
      \label{eqn:total-head-equation}
      \pdt u + \pdx \psi = - \frac{\qp{\fricR + \fricI +  k_0(u)}u}{h},
    \end{aligned}
  \end{equation}
  where $\psi$ is the total head, given by
  \begin{equation}
    \begin{aligned}
      \psi = \frac{u^2}{2} + \gravity h + \gravity Z.
    \end{aligned}
  \end{equation}
  Proceeding from this, multiplying the conservation of mass equation by $\psi$ we have
  \begin{equation}
    \begin{aligned}
      \psi\pdt h + \psi \pdx\qb{hu} = S\psi
    \end{aligned}
  \end{equation}
  which we can rewrite as
  \begin{equation}
    \begin{aligned}
      \qgroup{\pdt\qb{\psi h} + \pdx\qb{\psi hu}} - \qp{h\pdt\psi + \qp{hu}\pdx\psi} = S\psi.
    \end{aligned}
  \end{equation}
  The term $h\pdt\psi$ in the second component can be expanded as
  \begin{equation}
    \begin{aligned}
      h\pdt\psi
      &=
      \frac h2
      \pdt\qb{u^2}
      +
      h\gravity \pdt{ h}
      +
      h\gravity \pdt{ Z}
      \\
      &=
      \qp{hu} \pdt u
      +
      \frac\gravity2
      \pdt\qb{h^2}.
    \end{aligned}
  \end{equation}
  In the second component we may write
  \begin{equation}
    \begin{aligned}
      \pdx\psi = -\pdt u - \frac{\qp{\fricR + \fricI +  k_0(u)}u}{h}
    \end{aligned}
  \end{equation}
  from \eqref{eqn:total-head-equation}, and hence, after cancelling terms, we have
  \begin{equation}
    \begin{aligned}
      \qgroup{\pdt\qb{\psi h} + \pdx\qb{\psi hu}}
      -
      \qp{
        \pdt\qb{\gravity \frac{h^2}{2}}
        +
        \qp{\fricR + \fricI +  k_0(u)}
        u^2
      }
      =
      S\psi,
    \end{aligned}
  \end{equation}
  which can be rewritten as
  \begin{equation}
    \begin{aligned}
      \pdt\qb{\psi h - \gravity \frac{h^2}{2}} + \pdx\qb{\psi hu}
      =
      S\psi - \qp{\fricR + \fricI +  k_0(u)}u^2.
    \end{aligned}
  \end{equation}
  Expliciting $\psi$ as a function of $h,q,S,Z$ and oting that
  \begin{equation}
    \pdx\qb{\gravity\int_0^x S(s,t)Z(s)\d s}
    =
    \gravity S(t,x)Z(x)
    \Foreach (t,x)
  \end{equation}
  concludes the proof.
\end{Proofof}

\begin{Obs}[discontinuous solutions]
  In Theorems~\ref{THM_ESV} and \ref{THM_ESV2}, we only made reference
  to smooth solutions $(h, u)$ in defining the stability and entropy
  relations of the Saint-Venant system \eqref{ModelSL}.  For rough weak
  solutions, as with conservation laws, the entropy--entropy-flux pair
  $\qpreg{\hat E,\hat\Psi}$ has the potential to play a selection mechanism
  role to ensure uniqueness of weak solutions as it does with conservation
  laws.
\end{Obs}
\section{The numerical model}\label{SectionKinSch}
We now consider the numerical approximation of the Saint-Venant system
with rain. We follow the approach developed by
\citet{AudusseBristeauPerthame:00:techreport:Kinetic,PerthameSimeoni:01:article:A-kinetic},
with a suitable modification to accommodate the additional source and
friction terms. While any well-balanced computational finite volume method could be adapted
to simulate our model
\citep{Kroner:97:book:Numerical,LeVeque:92:book:Numerical,LeVeque:02:book:Finite,%
  Toro:09:book:Riemann,Bouchut:04:book:Nonlinear,Kurganov:18:article:Finite-volume},
the kinetic approach has the pleasant feature of naturally including
the additional term $Su$, and thereby also the
corresponding friction coefficients $\fricR$ and $\fricI$, in the
Saint-Venant system. As a result, as observed in
\citet{AudusseBristeauPerthame:00:techreport:Kinetic} the resulting schemes
are automatically up-winded and well balanced.
\subsection{Well balanced schemes}
A desirable property of the standard Saint-Venant system is the preservation of equilibrium states (referred to as the \emph{lake-at-rest}), given by
\begin{equation}
\begin{aligned}
h + Z = \const \tand u = 0.
\end{aligned}
\end{equation}
For our system, since through the addition of rainfall and infiltration effects we no longer have a conservation law but rather a balance law, we have the possibility of water being added to or lost from the lake, and thus this particular equilibrium only holds in the case $S = 0$. We must adapt this property therefore, and instead desire that our system preserves the \emph{filling-the-lake} state (see Figure~\ref{fig:filling-the-lake-state}), given by
\begin{equation}
\begin{aligned}
\pdt h = R \tand u = 0,
\end{aligned}
\end{equation}
that is, the rate at which the water height changes is equal to the rate at which water is added to the system through the rainfall term. Failing to preserve this property would mean a change is the mass of the water that is greater or lower than the rate at which it is added, thus violating the balance of mass property of our system.

\begin{figure}[ht]
\centering    
\scalebox{.8}{\begin{tikzpicture}
\begin{scope}
\clip (0,1.2) -- plot[smooth, domain=0:10](({\x}, {0.2*sin(deg(3*pi*\x/10))+1.5}) -- (10,4.5) -- (0,4.5);
\fill [cyan,opacity=0.3] (-0.5,-0.5) rectangle (10,5);
\end{scope}

\begin{scope}
\clip (0,0) -- plot[smooth, domain=0:10](({\x}, {0.2*sin(deg(3*pi*\x/10))+1.5}) -- (10,0);
\fill [brown,opacity=0.3] (-0.5,-0.5) rectangle (10,4);
\end{scope}

\draw [cyan, thick] plot[smooth, domain=0:10]({\x}, {4.5});
\draw [brown, thick] plot[smooth, domain=0:10]({\x}, {0.2*sin(deg(3*pi*\x/10))+1.5});
\draw [black, -latex] (-0.5,0) -- (10.5,0);
\draw [black, -latex] (0,0) -- (0,5.5);

\draw [black, latex-latex] (-0.25,0) -- node [left] {$Z(x)$} (-0.25,1.5);
\draw [black, latex-latex] (-0.25,1.5) -- node [left] {$h\xt$} (-0.25,4.5);
\draw [cyan] (5,2.9) node {wet region};
\draw [brown] (5,0.65) node {ground};


\draw [black, thick, -latex] (3.5,5.5) -- node [above left] {$R\xt$} (3.5,4.5);
\draw [black, thick, -latex] (6.5,5.5) -- node [above right] {$R\xt$} (6.5,4.5);

\draw (0,0) node [below] {0};
\draw (10,0) node [below] {L};
\end{tikzpicture}}
\caption{In considering a balance law system, we must adapt the \emph{lake-at-rest} property used for conservation laws to account for the addition of water.}
\label{fig:filling-the-lake-state}
\end{figure}

If we wish to maintain this property, we cannot rely on the usual
finite difference or finite volume methods, and thus a new (so called well-balanced)
scheme is required. Such an approach can be found by going back
to a kinetic interpretation of the system, as detailed in \cite{PerthameSimeoni:01:article:A-kinetic} and \cite{Ersoy:15:inproceedings:Dimension}. The method we use for the derivation of our kinetic scheme will follow much the same approach, though with the added complication of accounting for the additional terms. These \emph{kinetic solvers} can be modified to preserve the filling-the-lake state, while at the same time maintaining their simplicity and stability properties.

One of the direct benefits of using such an approach for the
Saint-Venant system is the ability of the kinetic solver to deal with
dry soil cases (that is, when $h = 0$)
\citet{AudusseBristeauPerthame:00:techreport:Kinetic}, which will be
of importance in ensuring our model continues to function if
infiltration causes the water level to fall close to zero, or if we
consider cases such as water flowing away from a beach front.
\subsection{Kinetic function}\label{subsec:kinetic-function}
To derive the kinetic equation for system \eqref{ModelSL}, we follow the kinetic formulation proposed by \citet{AudusseBristeauPerthame:00:techreport:Kinetic, PerthameSimeoni:01:article:A-kinetic} and further developed by
\cite{BourdariasErsoyGerbi:14:article:Unsteady, Ersoy:15:inproceedings:Dimension}. We consider a
\emph{kinetic averaging weight function} $\funk\chi\reals\reals$ and a
\emph{kinetic density function} $M$ satisfying
\begin{gather}
  \label{eq:kinetic-weight-properties}
  \chi(\omega) 
  = \chi(-\omega) \geq 0, 
  \quad 
  \int \chi(\omega)\d\omega 
  = 1, 
  \quad 
  \int \omega^2 \chi(\omega)\d\omega 
  = \frac{\gravity}{2}
  ,
  \\
  \label{eq:def:kinetic-density}
  M(t, x, \xi) 
  := 
  \sqrt{h(t, x)}\chi\qp{\frac{\xi - u(t, x)}{\sqrt{h(t, x)}}}.
\end{gather}
These functions originate in the kinetic theory where $M(t, x, \xi)$
accounts for the density of particles with speed $\xi$ at the
space-time point $(t, x)$.

In developing a numerical method, the goal
is for the derivation of the finite-volume scheme fluxes to be based
on $M$, through the following property which links the macroscopic
variables with the microscopic ones.

\begin{Pro}[macroscopic-microscopic relations]
\label{macroscopic-microscopic-relations}
Let the functions $h, u$ solve the Saint-Venant system
\eqref{ModelSL} and $M$ as in \eqref{eq:def:kinetic-density}. If $h(t,
x) > 0$ at $(t, x)$ then the following \emph{macroscopic-microscopic
  relations} hold
\begin{equation}
  \label{eq:macro-micro-relations}
  \int_{\reals} \discolvecthree[c]1\xi{\xi^2} M(t,x, \xi) \d\xi 
  = 
  \begin{bmatrix}[c] 
    h(t, x)
    \\ h(t, x)u(t, x)
    \\ 
    h(t, x)u(t, x)^2 + \fracl{\gravity h(t, x)^2}{2} 
  \end{bmatrix}
  .
\end{equation}
\end{Pro}

Recalling our Saint-Venant system \eqref{ModelSL}, we note that, by substituting $u = \sfrac{q}{h}$ (for $h > 0$), the topography and friction terms on the right-hand side of the conservation of momentum equation can be rewritten as
\begin{equation}
    -\gravity h\pdx  Z
    - \qp{\fricR + \fricI
      + k_0(u)}u
    \\
    =
    -\gravity h \qp{\pdx  Z
      + \frac{\qp{\fricR + \fricI + k_0(u)}u}{\gravity h }},
\end{equation}
following the approach considered in, for example,
\cite{BourdariasErsoyGerbi:11:article:A-kinetic}, adapted to account
for our additional friction terms. The reason for rewriting these
terms in this manner is so we can pack them into a single divergence
form; that is, we define the \indexemph{nonlinear flux integral
  operator}
\begin{equation}
  \label{eq:def:nonlinear-flux-integral-operator}
  \hat W[h(t,\cdot),u(t,\cdot)](x)
  :=
  Z(x)
  +
  \int_0^x \qb{\frac{\qp{\fricR + \fricI + k_0(u)}u}{\gravity h }}(t,s)
  \d s
\end{equation}
for each $x\in\reals$, and system \eqref{ModelSL} thus becomes
\begin{equation}
  \label{eq:Saint-Venant-with-rain-bar}
  \begin{aligned}
    \pdt  h
    +
    \pdx \qb{hu}
    &
    =
    S
    \\
    \pdt\qb{hu}
    +
    \pdx \qb{hu^2 + \frac{\gravity h^2}{2}}
    &
    =
    -\gravity h\pdx\hat W[h,u]
    + Su.
  \end{aligned}
\end{equation}
The kinetic scheme approach allows us to connect the Saint-Venant
system with the single scalar equation obtained by introducing an
auxilliary \indexemph{microscopic velocity} variable $\xi$
and looking at the evolution of the density
$\opinter0T\times\R2\ni(t,x,\xi)\mapsto M(t,x,\xi)$
as the solution of the following semilinear
\indexemph{kinetic equation}
\begin{equation}
  \label{eq:kinetic-equation}
  \pdt  M
  + \xi \pdx  M
  -
  \gravity \pdx\hat  W\qb{
    \average{M}_0,
    \frac{\average{M}_1}{\average{M}_0}
  }
  \partial_{\xi} M
  +
  \frac{
    S
    M
  }{
    \average{M}_0}
  = Q,
\end{equation}
where we use the following \indexemph{moment} notation for $m=0,1,\dotsc$
\begin{equation}
  \average{M}_m
  :=
  \int_{\reals} \xi^m
  M(\cdot,\cdot,\xi)\d\xi
  .
\end{equation}
The right-hand side in \eqref{eq:kinetic-equation},
$Q(t,x,\xi)$, plays the mathematical role of a
\emph{collision term}, similar, for instance, to the ones encountered
in Boltzmann's equation. In view of
Proposition~\ref{macroscopic-microscopic-relations}, if $(h,u)$
satisfying \eqref{ModelSL} is given, the pair $(M,Q)$ defined by
\eqref{eq:def:kinetic-density} and \eqref{eq:kinetic-equation} satisfy
the collision $0$-moment condition
\begin{equation}
  \label{eq:moment-condition}
  \average{Q}_m=0
  \text{ for }
  m=0,1.
\end{equation}
Conversely, each pair of functions $(M,Q)$ satisfying
\eqref{eq:kinetic-equation} and \eqref{eq:moment-condition} provides
a pair $(h,u)$ satisfying \eqref{ModelSL} by taking
\begin{equation}
  h:=\average{M}_0
  \tand
  hu:=\average{M}_1.
\end{equation}

We take a second here to note why the additional notation for $\average{M}_{0}$ and $\average{M}_{1}$ is required. In defining $\hat{W}$, we assume that $(h, u)$ are unknown and to be found, and thus we cannot use them within the definition. We are required, therefore, to instead define $\hat{W}$ in terms of $\average{M}_{0}$ and $\average{M}_{1}$, which are themselves defined in terms of the density function $M$, the solution to the kinetic equation, and which by virtue of the macroscopic-microscopic relations \eqref{eq:macro-micro-relations} will ultimately give us $h$ and $u$.

\begin{Obs}[advantages of the kinetic formulation]
  While the kinetic approach is one of many possibilities, and not
  without drawbacks, we metnion some of its nice features:
  \begin{enumerate}[label=(\roman*)]
  \item
    In contrast to previous work,
    e.g., \cite{PerthameSimeoni:01:article:A-kinetic}, the kinetic
    equation \eqref{eq:kinetic-equation} contains an extra term
    accounting for precipitation and infiltration effects. This
    departure is crucial for the derivation of the fluxes that lead to a
    well-balanced scheme in the presence of such terms.
  \item 
    We also note that, even though the Maxellian $M$ is constructed for still water
    steady states, where $S\xt = 0$, we can still use it here to ensure
    a well-balanced scheme.
  \item In general, it is easier to find a numerical scheme to solve
    equation \eqref{eq:kinetic-equation} for $M$ that has the
    properties we desire, such as entropy stability, than to solve the
    full Saint-Venant system for $h$ and $u$. However, in finding $M$,
    we can calculate $h$ and $hu$ by virtue of the macro-microscopic
    relations (proposition
    \ref{macroscopic-microscopic-relations}). In fact, $M$ is never
    calculated explicitly, rather the function
    \begin{equation}
      \begin{gathered}
        \label{eq:kinetic-Nemitskii}
        \hat M(\zeta,\varphi) := \sqrt \zeta \chi\qp{\frac{\varphi}{\sqrt\zeta}}\\
        \text{ whereby }
        M(t,x,\xi)=\hat M\qp{h(t,x),\xi-u(t,x)},
      \end{gathered}
    \end{equation}
    is used to build the fluxes appearing in a finite volume method, as we
    shall explain in \S\ref{sec:kinetic-flux}.
  \end{enumerate}
\end{Obs}
\subsection{Discretisation and kinetic fluxes}
\label{sec:kinetic-flux}
To go from the kinetic equation \eqref{eq:kinetic-equation} to a numerical method, we follow the approach of \cite{PerthameSimeoni:01:article:A-kinetic}, in which they developed a kinetic scheme for the standard Saint-Venant system, i.e. equation \eqref{eq:kinetic-equation} with $S = 0$. Their approach was based on the general method for developing a finite volume scheme, integrating the kinetic equation over the domain of interest, with the vector of unknowns defined as
\begin{equation}
  \begin{aligned}
    \label{eq:vector-unknowns}
    U_i^{n} = \int_{\reals} \begin{bmatrix} 1 \\ \xi \end{bmatrix} M_i^{n}(\xi) d\xi = \begin{bmatrix} h_i^n\\ h_i^n u_i^n \end{bmatrix},
  \end{aligned}
\end{equation}
where the final equality can be seen from the macroscopic-microscopic relations \eqref{eq:macro-micro-relations} and in view of the second observation above. We follow the same process for our Saint-Venant system, giving us the numerical scheme
\begin{equation}
  \label{eq:adapted-kinetic-scheme}
  \begin{aligned}
    U_i^{n+1} = U_i^n - \frac{\Delta t}{\Delta x}\qp{F_{i+1/2}^n - F_{i-1/2}^n} + \Delta t \begin{bmatrix} S_i^n \\ S_i^n u_i^n \end{bmatrix},
  \end{aligned}
\end{equation}
where $S_i^n$ is a discretisation of the combined rain and
infiltration terms. We pick the time-step, $\Delta t$,
according to
\begin{equation}
  \Delta t =
  \frac{\CFL\Delta x}{\max_i\qp{|u_i^n| + \sqrt{2 \gravity h_i^n}}},
\end{equation}
where $\CFL\in(0,1]$ is the Courant--Friedrichs--Levy stability constant
  \citep[e.g.]{PerthameSimeoni:01:article:A-kinetic}.
  
  The construction of the numerical fluxes $F_{i\pm1/2}^n$ in
  \eqref{eq:adapted-kinetic-scheme} is based on the operator
  associated with $\hat M$ given in \eqref{eq:kinetic-Nemitskii}. We
  give here a brief overview of how the successive numerical flux
  terms are developed; the interested reader is directed to
  \cite{PerthameSimeoni:01:article:A-kinetic} and
  \cite{BourdariasErsoyGerbi:14:article:Unsteady} for more details.
  
  We define the numerical flux (which has two components representing the conservation of mass and momentum equations), as the integral
\begin{equation}
\label{eq:numerical-flux-integral}
F_{i\pm\frac12}^n := \int_\reals \xi \begin{bmatrix}1 \\ \xi \end{bmatrix} M_{i\pm1/2}^\mp(\xi) d\xi.
\end{equation}
The intermediate quantities $M_{i\pm1/2}^{\mp}(\xi)$, which measure the flux at the upper and lower boundaries of the cell $c_i = [x_{i-1/2}, x_{i+1/2}]$, respectively, are realised as upwinded fluxes:
\begin{equation}
  \begin{split}
    \label{eq:numerical-flux-defn}
    M_{i+1/2}^- &:= M_i^n(\xi)\iverson {\xi>0} + M_{i+1/2}^n(\xi)\iverson{\xi<0},
    \\
    M_{i-1/2}^+ &:= M_i^n(\xi) \iverson {\xi<0} + M_{i-1/2}^n(\xi)\iverson{\xi>0},
  \end{split}
\end{equation}
with
\begin{equation}
  \begin{aligned}
    \label{eq:numerical-flux-defn-intermediate}
    M_{i\pm1/2}^n &:= M_i^n(-\xi) \iverson{\abs{\xi}^2\leq2\gravity \Delta W_{i\pm1/2}^n}
    \\
    &\qquad + M_{i\pm1}^n \qp{\mp{\sqrt{\abs{\xi}^2-2\gravity \Delta W_{i\pm1/2}^n}}} \iverson{\abs{\xi}^2\geq 2\gravity \Delta W_{i\pm1/2}^n},
  \end{aligned}
\end{equation}
where we use the notation
\begin{equation}
\begin{aligned}
  \iverson P
  :=
  \begin{cases}
    1 & \text{if P is true},
    \\
    0 & \text{if P is false}.
  \end{cases}
\end{aligned}
\end{equation}
To understand how these fluxes are defined, consider the flux over the
upper bound $M_{i+1/2}^-$. From \eqref{eq:numerical-flux-defn}, we can
see that it comprises two terms:
\begin{enumerate}[label=(\roman*)]
\item
  $M_i^n(\xi)\iverson {\xi>0}$: movement of water with positive
  velocity ($\xi>0$) from within cell $c_i$ to cell $c_{i+1}$
\item
  $M_{i+1/2}^n(\xi)\iverson{\xi<0}$: movement of water with negative
  velocity ($\xi<0$) from within cell $c_{i+1}$ to cell $c_{i}$. This
  term is decomposed a second time into components reflecting whether
  the water has enough energy to overcome the topography and friction
  to enter or leave the cell.
\end{enumerate}
A similar decomposition exists for the second term $M_{i-1/2}^+$, where this time we consider the negative and positive velocity of $\xi$ for each case, respectively.

The term $\Delta W_{i\pm1/2}^n$ is the upwinded source term and provides the jump condition necessary for a particle in one cell to overcome the friction and topography to move to an adjacent cell. Consistent with previous definitions, we calculate this term numerically as:
\begin{equation}
\begin{aligned}
\label{eq:numerical-packed-source-term}
\Delta W_{i+1/2}^n = W_{i+1}(t_n) - W_i(t_n), \quad \Delta W_{i-1/2}^n = W_{i-1}(t_n) - W_i(t_n),
\end{aligned}
\end{equation}
where
\begin{equation}
\begin{aligned}
W_i(t) = \iverson{c_i}(x)\frac{1}{\Delta x}\int_{c_i}W\xt dx
\end{aligned}
\end{equation}
for a given cell $c_i$. The semi-discretised kinetic density, $M_i^n$, is defined as
\begin{equation}
\begin{aligned}
M_i^n(\xi) := \sqrt{h_i^n} \, \chi\qp{\frac{\xi - u_i^n}{\sqrt{h_i^n}}}.
\end{aligned}
\end{equation}

The discretisation we use in our scheme will be based upon the \indexemph{Barrenblatt kinetic weighting function}
\begin{equation}
  \chi(\omega) = \frac{1}{\pi \gravity}\sqrt{\qp{2\gravity - \omega^2}_+}
  \text{ for }
  \omega\in\reals,
\end{equation}
where $(X)_{+}$ stands for the positive part of $X$ \citep[eq.(2.13)]{PerthameSimeoni:01:article:A-kinetic}. We note that this choice of function satisfies the properties we outlined in \S\ref{subsec:kinetic-function}.

\section{Numerical tests}
\label{SectionNum}
The kinetic scheme we use for our numerical method was implemented by extending the code of \cite{BessonLakkis:13:url:Finite} to account for the additional source term in \eqref{eq:adapted-kinetic-scheme}, and we present here two simple numerical tests to demonstrate the validity and application of our Saint-Venant system and the associated numerical method. For this paper, we will focus on the rain term only and thus assume $I \equiv 0$; the coupling of a realistic infiltration model with our Saint-Venant system and the treatment of the boundary conditions would be a paper in itself, and thus will be considered in future research.

  In \S\ref{subsec:real-world-experiment}, we compare the accuracy of our numerical model to a flume experiment, considering how our results compare to both the physical data collected from this experiment and also previous numerical simulations of this experiment by other authors. Then, in \S\ref{subsec:cascade}, we simulate multiple rainfall processes of increasing duration on a slope with both a constant and decreasing gradient, and measure how the value of the friction effect $\alpha$ impacts the solution.

  \subsection{Comparison with real-world data}
  \label{subsec:real-world-experiment}
  For our first test, we explore the accuracy of our numerical scheme by comparing with data taken from the flume experiment run as part of the ANR project METHODE at INRA-Orl\'eans; we also compare our results to those obtained by \citet{DelestreJames:08:inproceedings:Simulation}, who considered only an addition to the conservation of mass equation.
  \begin{figure}[ht]
    \centering   
    \scalebox{1.0}{\begin{tikzpicture}
        \draw [black] (0,0) -- (8,0);
        \draw [black, latex-latex] (0,-0.1) -- node [below] {$4$m} (8,-0.1);

        \draw [black] (0,0) -- (0,2);
        \draw [black, latex-latex] (-0.1,0) -- node [left] {$0.2$m} (-0.1,2);
        \draw [black] (0,2) -- (8,0);

        \begin{scope}
          \clip (0,2) -- (0,2.1) -- (8,0.1) -- (8,0) -- (0,2);
          \fill [cyan,opacity=0.3] (-0.5,-0.5) rectangle (8,2.1);
        \end{scope}

        \draw [cyan] (0,2) -- (0,2.1) -- (8,0.1) -- (8,0);
        \draw [black] (0,3) -- (7.5,3);
        \foreach \x in {0, 0.5, 1, ..., 7.5} {\draw [cyan, thick, dashed, -latex] (\x, 3) -- (\x, 2.5);}
    \end{tikzpicture}}
    \caption{Visualisation of the flume experiment.}
  \end{figure}

  The experiment in question concerns a slope with a $5\%$ gradient, an initial height $h_0 = 0$, and initial discharge $q_0 = 0$. The topography for the slope is given by
  \begin{equation}
    \begin{aligned}
      Z(x) &= 0.2 - \frac{x}{20}\quad \forall x \in [0,4],
    \end{aligned}
  \end{equation}
  we consider $N = 1000$ meshpoints, and we assume a CFL number of $0.95$. Rain falls onto the slope uniformly at a constant rate within a given time interval,
  \begin{equation}
    \begin{aligned}
      R(t,x) =
      \begin{cases}
        50 \text{ mm/hr} & \text{if } (t, x) \in [5,125] \times [0,3.95],\\
        0 & \text{otherwise},
      \end{cases}
    \end{aligned}
  \end{equation}
  and we measure the discharge at the downstream edge of the slope up to time $T = 250$s. For our simulation, we assume the rain-induced friction level $\alpha = 1$. The hydrograph for the experiment, together with the simulated values, is provided below in Figure~\ref{fig:hydrograph-flume-experiment}.
  \ThesisSimplify{\boolean{simplify}}{Flume experiment}{
    \begin{figure}[h]
      \centering
      \scalebox{1.0}{\begin{tikzpicture}
          \begin{axis}[width=1.0\textwidth, height=0.55\textwidth, xlabel=time, ylabel=discharge $(g/s)$, legend pos=north east, legend cell align={left}, ymin=0, ymax=9, xmin=0, xmax=250]
            \addplot [blue, fill=blue, mark=x, mark size=2] table [scatter, only marks, x=time, y=median, col sep=comma] {./Experiment_full.txt};
            \addplot [red, thick] table [x=time, y=discharge, col sep=comma] {./Simulation_half.txt};
            \addlegendentry{\ experiment}
            \addlegendentry{\ simulation}
          \end{axis}
      \end{tikzpicture}}
      \caption{Hydrograph for the uniform slope test with both the experimental data and the simulated solution.}
      \label{fig:hydrograph-flume-experiment}
  \end{figure}}

  The results from the simulation compare well with with the experimental data, particularly in the third and final phase when no rain is falling and the discharge is decreasing gradually over time. For the initial phase, where the discharge is increasing, the simulation matches well to begin with (up to around $25$s) but subsequently appears to increase at a slower rate than in the experiment; at $t = 35$s, for example, the experiment shows a discharge level of $5.7$ g/s, compared to the simulation which is at $3.59$ g/s. For the secondary phase, where the discharge has stabilised, while the simulation does not capture the fluctuations that were present in the experiment, it does maintain a smooth transition through the centre of the data, and begins to decrease at the same point as the experiment.

  Comparing the results of our simulation to those of \citet{DelestreJames:08:inproceedings:Simulation}, we see that our simulation matches much better in both the initial and final phases; in the initial phase, their simulation increases a lot sooner than in the experiment, and though it begins to decrease at the correct time, it also underestimates the total discharge level in the final phase. For the secondary phase, their results are slightly lower than ours but still consistent with the experimental data.

  \subsection{Single-level and three-level cascade}
  \label{subsec:cascade}
  For our second test, we consider a similar scenario to \S\ref{subsec:real-world-experiment} but with a higher intensity rainfall-runoff process on a slope with a much shallower gradient. We compare how the water flows when the gradient of the slope is constant, and how it flows when the gradient decreases periodically from the upstream to downstream end. We consider a spatial domain $x \in [0, 12]$, a final time $t = 40$s, a rainfall intensity $R_0 = 0.001$ which falls across the entire domain, $N = 1000$ meshpoints, and a CFL number of $0.95$. The parameters that we want to change and measure the effect of are the following:
  \begin{enumerate}[label=(\arabic*)]
  \item the total length of the rainfall process $T_R = 10, 20$, and $30$ seconds
  \item the topography of the slope onto which the rain falls, for which we consider a constant slope (the single cascade) with $Z_1(x) = (12 - x)0.005$, and a decreasing slope (the three-level cascade, see Figure~\ref{fig:three-level-cascade}) with
    \begin{equation}
      \begin{aligned}
        Z_2(x) &= \begin{cases} (12 - x)0.006 - 0.012 & \text{if } x \in [0,4] \\
	  (12 - x)0.005 - 0.004 & \text{if } x \in [4,8] \\
	  (12 - x)0.004 & \text{if } x \in [8,12] \end{cases}
      \end{aligned}
    \end{equation}
  \item the rain-induced friction level $\alpha$, for which we take $\alpha = 0, 1,$ and $5$ for the single cascade and $\alpha = 0$ and $1$ for the three-level cascade
  \end{enumerate}
  \begin{figure}[ht]
    \centering    
    \scalebox{0.75}{\begin{tikzpicture}
        \draw [black] (0,0) -- (12,0);
        \draw [black, latex-latex] (0,-0.1) -- node [below] {$4$m} (4,-0.1);
        \draw [black, latex-latex] (4,-0.1) -- node [below] {$4$m} (8,-0.1);
        \draw [black, latex-latex] (8,-0.1) -- node [below] {$4$m} (12,-0.1);

        \draw [black] (12,0) -- (8,0.4) -- (4,1.2) -- (0,2.8);

        \begin{scope}
          \clip (12,0) -- (12,0.1) -- (8,0.5) -- (4,1.3) -- (0,2.9) -- (0,2.8) -- (4,1.2) -- (8,0.4) -- (12,0);
          \fill [cyan,opacity=0.3] (-0.5,-0.5) rectangle (12,5.7);
        \end{scope}

        \draw [cyan] (12,0) -- (12,0.1) -- (8,0.5) -- (4,1.3) -- (0,2.9) -- (0,2.8);
        \draw [black, dashed] (0,0) -- (0,2.8);
        \draw [black, dashed] (4,0) -- (4,1.2);
        \draw [black, dashed] (8,0) -- (8,0.4);

        \draw [black] (0,5) -- (12,5);
        \foreach \x in {0, 0.5, ..., 12} {\draw [cyan, thick, dashed, -latex] (\x, 5) -- (\x, 4);}
    \end{tikzpicture}}
    \caption{Topography for the three-level cascade, showing the decrease in gradient from the upstream to downstream end.}
    \label{fig:three-level-cascade}
  \end{figure}

  For our outputs, we measure the height of the flow across the entire domain at the point the rainfall stops, and we also measure the height and discharge at the downstream end (i.e. $x = 12$) up to the final time.

  Starting with the height profile at $t = T_R$, we see in Figure~\ref{fig:single-height} that for the single-level cascade, increasing the value of $\alpha$ causes more water to accumulate at the upstream end, as expected since the momentum of the water will decrease the higher the value of $\alpha$. It is interesting to note also that for $T_R = 10$, the water still accumulates to the same maximum amount irrespective of the value of $\alpha$, though at different points; for longer rainfall times, this accumulation can still occur, but potentially beyond the end of the domain if the water has enough momentum.

  \ThesisSimplify{\boolean{simplify}}{Height single-cascade}{
    \begin{figure}[h]
      \centering
      \scalebox{0.55}{\begin{tikzpicture}
          \begin{axis}[width=0.8\textwidth, height=0.5\textwidth, xlabel=x, ylabel=height $h$, legend pos=north west, legend cell align={left}, ymin=0, ymax=0.035, xmin=0, xmax=12, yticklabel style={/pgf/number format/.cd, fixed, fixed zerofill, precision=2, /tikz/.cd}, scaled y ticks = false]
            \addplot [red, thick] table [x=time, y=height, col sep=comma] {./PaperGraphs/Test3_time10_height_a0_single.txt};
            \addplot [blue, thick] table [x=time, y=height, col sep=comma] {./PaperGraphs/Test3_time10_height_a1_single.txt};
            \addplot [green, thick] table [x=time, y=height, col sep=comma] {./PaperGraphs/Test3_time10_height_a5_single.txt};
            \addlegendentry{\ $\alpha = 0$}
            \addlegendentry{\ $\alpha = 1$}
            \addlegendentry{\ $\alpha = 5$}
          \end{axis}
      \end{tikzpicture}}
      \scalebox{0.55}{\begin{tikzpicture}
          \begin{axis}[width=0.8\textwidth, height=0.5\textwidth, xlabel=x, ylabel=height $h$, legend pos=north west, legend cell align={left}, ymin=0, ymax=0.035, xmin=0, xmax=12, yticklabel style={/pgf/number format/.cd, fixed, fixed zerofill, precision=2, /tikz/.cd}, scaled y ticks = false]
            \addplot [red, thick] table [x=time, y=height, col sep=comma] {./PaperGraphs/Test3_time20_height_a0_single.txt};
            \addplot [blue, thick] table [x=time, y=height, col sep=comma] {./PaperGraphs/Test3_time20_height_a1_single.txt};
            \addplot [green, thick] table [x=time, y=height, col sep=comma] {./PaperGraphs/Test3_time20_height_a5_single.txt};
            \addlegendentry{\ $\alpha = 0$}
            \addlegendentry{\ $\alpha = 1$}
            \addlegendentry{\ $\alpha = 5$}
          \end{axis}
      \end{tikzpicture}}
      \scalebox{0.55}{\begin{tikzpicture}
          \begin{axis}[width=0.8\textwidth, height=0.5\textwidth, xlabel=x, ylabel=height $h$, legend pos=north west, legend cell align={left}, ymin=0, ymax=0.035, xmin=0, xmax=12, yticklabel style={/pgf/number format/.cd, fixed, fixed zerofill, precision=2, /tikz/.cd}, scaled y ticks = false]
            \addplot [red, thick] table [x=time, y=height, col sep=comma] {./PaperGraphs/Test3_time30_height_a0_single.txt};
            \addplot [blue, thick] table [x=time, y=height, col sep=comma] {./PaperGraphs/Test3_time30_height_a1_single.txt};
            \addplot [green, thick] table [x=time, y=height, col sep=comma] {./PaperGraphs/Test3_time30_height_a5_single.txt};
            \addlegendentry{\ $\alpha = 0$}
            \addlegendentry{\ $\alpha = 1$}
            \addlegendentry{\ $\alpha = 5$}
          \end{axis}
      \end{tikzpicture}}
      \caption{Height profiles for the single-level cascade for varying values of $\alpha$ at the rainfall end time $T_R = 10$, $20$, and $30$.}
      \label{fig:single-height}
  \end{figure}}
  
  For the three-level cascade shown in Figure~\ref{fig:triple-height}, the reduction in gradient induces multiple waves to be formed, though these waves become more smoothed out as the rainfall time increases. This effect is reduced as $\alpha$ is increased, consistent with our expectation since the water flow will be more slowed. The flows do not accumulate to the same amount, though this may be due to the length of the domain.
  \ThesisSimplify{\boolean{simplify}}{Height triple-cascade}{
    \begin{figure}[h]
      \centering
      \scalebox{0.55}{\begin{tikzpicture}
          \begin{axis}[width=0.8\textwidth, height=0.5\textwidth, xlabel=x, ylabel=height $h$, legend pos=north west, legend cell align={left}, ymin=0, ymax=0.035, xmin=0, xmax=12, yticklabel style={/pgf/number format/.cd, fixed, fixed zerofill, precision=2, /tikz/.cd}, scaled y ticks = false]
            \addplot [red, thick] table [x=time, y=height, col sep=comma] {./PaperGraphs/Test3_time10_height_a0_triple.txt};
            \addplot [blue, thick] table [x=time, y=height, col sep=comma] {./PaperGraphs/Test3_time10_height_a1_triple.txt};
            \addlegendentry{\ $\alpha = 0$}
            \addlegendentry{\ $\alpha = 1$}
          \end{axis}
      \end{tikzpicture}}
      \scalebox{0.55}{\begin{tikzpicture}
          \begin{axis}[width=0.8\textwidth, height=0.5\textwidth, xlabel=x, ylabel=height $h$, legend pos=north west, legend cell align={left}, ymin=0, ymax=0.035, xmin=0, xmax=12, yticklabel style={/pgf/number format/.cd, fixed, fixed zerofill, precision=2, /tikz/.cd}, scaled y ticks = false]
            \addplot [red, thick] table [x=time, y=height, col sep=comma] {./PaperGraphs/Test3_time20_height_a0_triple.txt};
            \addplot [blue, thick] table [x=time, y=height, col sep=comma] {./PaperGraphs/Test3_time20_height_a1_triple.txt};
            \addlegendentry{\ $\alpha = 0$}
            \addlegendentry{\ $\alpha = 1$}
          \end{axis}
      \end{tikzpicture}}
      \scalebox{0.55}{\begin{tikzpicture}
          \begin{axis}[width=0.8\textwidth, height=0.5\textwidth, xlabel=x, ylabel=height $h$, legend pos=north west, legend cell align={left}, ymin=0, ymax=0.035, xmin=0, xmax=12, yticklabel style={/pgf/number format/.cd, fixed, fixed zerofill, precision=2, /tikz/.cd}, scaled y ticks = false]
            \addplot [red, thick] table [x=time, y=height, col sep=comma] {./PaperGraphs/Test3_time30_height_a0_triple.txt};
            \addplot [blue, thick] table [x=time, y=height, col sep=comma] {./PaperGraphs/Test3_time30_height_a1_triple.txt};
            \addlegendentry{\ $\alpha = 0$}
            \addlegendentry{\ $\alpha = 1$}
          \end{axis}
      \end{tikzpicture}}
      \caption{Height profiles for the three-level cascade for varying values of $\alpha$ at the rainfall end time $T_R = 10$, $20$, and $30$.}
      \label{fig:triple-height}
  \end{figure}}

  For the second part of our numerical test, we consider the height and momentum profiles over time at the end of the domain ($x = 12$). For the single-level cascade (see Figure~\ref{fig:single-height-momentum}), we see that increasing the friction level $\alpha$ extends the height profile of the flow, causing it to decrease at a later time. This occurs for all lengths of rainfall $T_R$, though notably we see that as $T_R$ increases the length of time for which the flow plateaus is decreased.

  For the momentum, the change in friction and length of rainfall has a much more pronounced effect. For $T_R = 10$, the three friction levels result in much the same profile though slightly shifted as the level increases. For $T_R = 20$ and $30$, however, the profiles are more varied, with the momentum tending to change rather linearly for $\alpha = 0$ but showing a much more curved profile for $\alpha = 5$. We can also see that for $T_R = 30$, the momentum is decreasing when the rainfall stops for $\alpha = 0$ and $1$, but continues to rise and at an increased rate for $\alpha = 5$.
  \ThesisSimplify{\boolean{simplify}}{Height and momentum single-cascade}{
    \begin{figure}[h]
      \centering
      \scalebox{0.6}{\begin{tikzpicture}
          \begin{axis}[width=0.775\textwidth, height=0.5\textwidth, xlabel=t, ylabel=height $h$, title={$T_R = 10$}, legend pos=north west, legend cell align={left}, ymin=0, ymax=0.02, xmin=0, xmax=40, yticklabel style={/pgf/number format/.cd, fixed, fixed zerofill, precision=3, /tikz/.cd}, scaled y ticks = false, y label style={at={(axis description cs:-0.025,0.5)}}]
            \addplot [red, thick] table [x=time, y=height, col sep=comma] {./PaperGraphs/Test3_time10_discharge_a0_single.txt};
            \addplot [blue, thick] table [x=time, y=height, col sep=comma] {./PaperGraphs/Test3_time10_discharge_a1_single.txt};
            \addplot [green, thick] table [x=time, y=height, col sep=comma] {./PaperGraphs/Test3_time10_discharge_a5_single.txt};
            \addlegendentry{\ $\alpha = 0$}
            \addlegendentry{\ $\alpha = 1$}
            \addlegendentry{\ $\alpha = 5$}
          \end{axis}
      \end{tikzpicture}}
      \scalebox{0.6}{\begin{tikzpicture}
          \begin{axis}[width=0.775\textwidth, height=0.5\textwidth, xlabel=t, ylabel=momentum $q$, title={$T_R = 10$}, legend pos=north west, legend cell align={left}, ymin=0, ymax=0.01, xmin=0, xmax=40, yticklabel style={/pgf/number format/.cd, fixed, fixed zerofill, precision=3, /tikz/.cd}, scaled y ticks = false, y label style={at={(axis description cs:-0.025,0.5)}}]
            \addplot [red, thick] table [x=time, y=discharge, col sep=comma] {./PaperGraphs/Test3_time10_discharge_a0_single.txt};
            \addplot [blue, thick] table [x=time, y=discharge, col sep=comma] {./PaperGraphs/Test3_time10_discharge_a1_single.txt};
            \addplot [green, thick] table [x=time, y=discharge, col sep=comma] {./PaperGraphs/Test3_time10_discharge_a5_single.txt};
            \addlegendentry{\ $\alpha = 0$}
            \addlegendentry{\ $\alpha = 1$}
            \addlegendentry{\ $\alpha = 5$}
          \end{axis}
      \end{tikzpicture}}
      \scalebox{0.6}{\begin{tikzpicture}
          \begin{axis}[width=0.775\textwidth, height=0.5\textwidth, xlabel=t, ylabel=height $h$, title={$T_R = 20$}, legend pos=north west, legend cell align={left}, ymin=0, ymax=0.025, xmin=0, xmax=40, yticklabel style={/pgf/number format/.cd, fixed, fixed zerofill, precision=3, /tikz/.cd}, scaled y ticks = false, y label style={at={(axis description cs:-0.025,0.5)}}]
            \addplot [red, thick] table [x=time, y=height, col sep=comma] {./PaperGraphs/Test3_time20_discharge_a0_single.txt};
            \addplot [blue, thick] table [x=time, y=height, col sep=comma] {./PaperGraphs/Test3_time20_discharge_a1_single.txt};
            \addplot [green, thick] table [x=time, y=height, col sep=comma] {./PaperGraphs/Test3_time20_discharge_a5_single.txt};
            \addlegendentry{\ $\alpha = 0$}
            \addlegendentry{\ $\alpha = 1$}
            \addlegendentry{\ $\alpha = 5$}
          \end{axis}
      \end{tikzpicture}}
      \scalebox{0.6}{\begin{tikzpicture}
          \begin{axis}[width=0.775\textwidth, height=0.5\textwidth, xlabel=t, ylabel=momentum $q$, title={$T_R = 20$}, legend pos=north west, legend cell align={left}, ymin=0, ymax=0.02, xmin=0, xmax=40, yticklabel style={/pgf/number format/.cd, fixed, fixed zerofill, precision=3, /tikz/.cd}, scaled y ticks = false, y label style={at={(axis description cs:-0.025,0.5)}}]
            \addplot [red, thick] table [x=time, y=discharge, col sep=comma] {./PaperGraphs/Test3_time20_discharge_a0_single.txt};
            \addplot [blue, thick] table [x=time, y=discharge, col sep=comma] {./PaperGraphs/Test3_time20_discharge_a1_single.txt};
            \addplot [green, thick] table [x=time, y=discharge, col sep=comma] {./PaperGraphs/Test3_time20_discharge_a5_single.txt};
            \addlegendentry{\ $\alpha = 0$}
            \addlegendentry{\ $\alpha = 1$}
            \addlegendentry{\ $\alpha = 5$}
          \end{axis}
      \end{tikzpicture}}
      \scalebox{0.6}{\begin{tikzpicture}
          \begin{axis}[width=0.775\textwidth, height=0.5\textwidth, xlabel=t, ylabel=height $h$, title={$T_R = 30$}, legend pos=north west, legend cell align={left}, ymin=0, ymax=0.035, xmin=0, xmax=40, yticklabel style={/pgf/number format/.cd, fixed, fixed zerofill, precision=3, /tikz/.cd}, scaled y ticks = false, y label style={at={(axis description cs:-0.025,0.5)}}]
            \addplot [red, thick] table [x=time, y=height, col sep=comma] {./PaperGraphs/Test3_time30_discharge_a0_single.txt};
            \addplot [blue, thick] table [x=time, y=height, col sep=comma] {./PaperGraphs/Test3_time30_discharge_a1_single.txt};
            \addplot [green, thick] table [x=time, y=height, col sep=comma] {./PaperGraphs/Test3_time30_discharge_a5_single.txt};
            \addlegendentry{\ $\alpha = 0$}
            \addlegendentry{\ $\alpha = 1$}
            \addlegendentry{\ $\alpha = 5$}
          \end{axis}
      \end{tikzpicture}}
      \scalebox{0.6}{\begin{tikzpicture}
          \begin{axis}[width=0.775\textwidth, height=0.5\textwidth, xlabel=t, ylabel=momentum $q$, title={$T_R = 30$}, legend pos=north west, legend cell align={left}, ymin=0, ymax=0.02, xmin=0, xmax=40, yticklabel style={/pgf/number format/.cd, fixed, fixed zerofill, precision=3, /tikz/.cd}, scaled y ticks = false, y label style={at={(axis description cs:-0.025,0.5)}}]
            \addplot [red, thick] table [x=time, y=discharge, col sep=comma] {./PaperGraphs/Test3_time30_discharge_a0_single.txt};
            \addplot [blue, thick] table [x=time, y=discharge, col sep=comma] {./PaperGraphs/Test3_time30_discharge_a1_single.txt};
            \addplot [green, thick] table [x=time, y=discharge, col sep=comma] {./PaperGraphs/Test3_time30_discharge_a5_single.txt};
            \addlegendentry{\ $\alpha = 0$}
            \addlegendentry{\ $\alpha = 1$}
            \addlegendentry{\ $\alpha = 5$}
          \end{axis}
      \end{tikzpicture}}
      \caption{Height and momentum profiles over time for the single-level cascade for varying values of rain-induced friction $\alpha$ at the domain end $x = 12$.}
      \label{fig:single-height-momentum}
  \end{figure}}

  The three-level cascade shown in Figure~\ref{fig:triple-height-momentum} exhibits much the same behaviour as the single-level cascade, with perhaps the most notable changes being in the height profile; for the single-level cascade the height profile for $\alpha = 0$ remains consistently below that for $\alpha = 1$, but for the three-level cascade this does not hold true between, approximately, $t = 10$ and $t = 20$. The momentum profile for the three-level cascade is very similar in behaviour to the single-level.
  \ThesisSimplify{\boolean{simplify}}{Height and momentum triple-cascade}{
    \begin{figure}[h]
      \centering
      \scalebox{0.6}{\begin{tikzpicture}
          \begin{axis}[width=0.775\textwidth, height=0.5\textwidth, xlabel=t, ylabel=height $h$, title={$T_R = 10$}, legend pos=north west, legend cell align={left}, ymin=0, ymax=0.02, xmin=0, xmax=40, yticklabel style={/pgf/number format/.cd, fixed, fixed zerofill, precision=3, /tikz/.cd}, scaled y ticks = false, y label style={at={(axis description cs:-0.025,0.5)}}]
            \addplot [red, thick] table [x=time, y=height, col sep=comma] {./PaperGraphs/Test3_time10_discharge_a0_triple.txt};
            \addplot [blue, thick] table [x=time, y=height, col sep=comma] {./PaperGraphs/Test3_time10_discharge_a1_triple.txt};
            \addlegendentry{\ $\alpha = 0$}
            \addlegendentry{\ $\alpha = 1$}
          \end{axis}
      \end{tikzpicture}}
      \scalebox{0.6}{\begin{tikzpicture}
          \begin{axis}[width=0.775\textwidth, height=0.5\textwidth, xlabel=t, ylabel=momentum $q$, title={$T_R = 10$}, legend pos=north west, legend cell align={left}, ymin=0, ymax=0.015, xmin=0, xmax=40, yticklabel style={/pgf/number format/.cd, fixed, fixed zerofill, precision=3, /tikz/.cd}, scaled y ticks = false, y label style={at={(axis description cs:-0.025,0.5)}}]
            \addplot [red, thick] table [x=time, y=discharge, col sep=comma] {./PaperGraphs/Test3_time10_discharge_a0_triple.txt};
            \addplot [blue, thick] table [x=time, y=discharge, col sep=comma] {./PaperGraphs/Test3_time10_discharge_a1_triple.txt};
            \addlegendentry{\ $\alpha = 0$}
            \addlegendentry{\ $\alpha = 1$}
          \end{axis}
      \end{tikzpicture}}
      \scalebox{0.6}{\begin{tikzpicture}
          \begin{axis}[width=0.775\textwidth, height=0.5\textwidth, xlabel=t, ylabel=height $h$, title={$T_R = 20$}, legend pos=north west, legend cell align={left}, ymin=0, ymax=0.03, xmin=0, xmax=40, yticklabel style={/pgf/number format/.cd, fixed, fixed zerofill, precision=3, /tikz/.cd}, scaled y ticks = false, y label style={at={(axis description cs:-0.025,0.5)}}]
            \addplot [red, thick] table [x=time, y=height, col sep=comma] {./PaperGraphs/Test3_time20_discharge_a0_triple.txt};
            \addplot [blue, thick] table [x=time, y=height, col sep=comma] {./PaperGraphs/Test3_time20_discharge_a1_triple.txt};
            \addlegendentry{\ $\alpha = 0$}
            \addlegendentry{\ $\alpha = 1$}
          \end{axis}
      \end{tikzpicture}}
      \scalebox{0.6}{\begin{tikzpicture}
          \begin{axis}[width=0.775\textwidth, height=0.5\textwidth, xlabel=t, ylabel=momentum $q$, title={$T_R = 20$}, legend pos=north west, legend cell align={left}, ymin=0, ymax=0.02, xmin=0, xmax=40, yticklabel style={/pgf/number format/.cd, fixed, fixed zerofill, precision=3, /tikz/.cd}, scaled y ticks = false, y label style={at={(axis description cs:-0.025,0.5)}}]
            \addplot [red, thick] table [x=time, y=discharge, col sep=comma] {./PaperGraphs/Test3_time20_discharge_a0_triple.txt};
            \addplot [blue, thick] table [x=time, y=discharge, col sep=comma] {./PaperGraphs/Test3_time20_discharge_a1_triple.txt};
            \addlegendentry{\ $\alpha = 0$}
            \addlegendentry{\ $\alpha = 1$}
          \end{axis}
      \end{tikzpicture}}
      \scalebox{0.6}{\begin{tikzpicture}
          \begin{axis}[width=0.775\textwidth, height=0.5\textwidth, xlabel=t, ylabel=height $h$, title={$T_R = 30$}, legend pos=north west, legend cell align={left}, ymin=0, ymax=0.03, xmin=0, xmax=40, yticklabel style={/pgf/number format/.cd, fixed, fixed zerofill, precision=3, /tikz/.cd}, scaled y ticks = false, y label style={at={(axis description cs:-0.025,0.5)}}]
            \addplot [red, thick] table [x=time, y=height, col sep=comma] {./PaperGraphs/Test3_time30_discharge_a0_triple.txt};
            \addplot [blue, thick] table [x=time, y=height, col sep=comma] {./PaperGraphs/Test3_time30_discharge_a1_triple.txt};
            \addlegendentry{\ $\alpha = 0$}
            \addlegendentry{\ $\alpha = 1$}
          \end{axis}
      \end{tikzpicture}}
      \scalebox{0.6}{\begin{tikzpicture}
          \begin{axis}[width=0.775\textwidth, height=0.5\textwidth, xlabel=t, ylabel=momentum $q$, title={$T_R = 30$}, legend pos=north west, legend cell align={left}, ymin=0, ymax=0.02, xmin=0, xmax=40, yticklabel style={/pgf/number format/.cd, fixed, fixed zerofill, precision=3, /tikz/.cd}, scaled y ticks = false, y label style={at={(axis description cs:-0.025,0.5)}}]
            \addplot [red, thick] table [x=time, y=discharge, col sep=comma] {./PaperGraphs/Test3_time30_discharge_a0_triple.txt};
            \addplot [blue, thick] table [x=time, y=discharge, col sep=comma] {./PaperGraphs/Test3_time30_discharge_a1_triple.txt};
            \addlegendentry{\ $\alpha = 0$}
            \addlegendentry{\ $\alpha = 1$}
          \end{axis}
      \end{tikzpicture}}
      \caption{Height and momentum profiles over time for the three-level cascade for varying values of rain-induced friction $\alpha$ at the domain end $x = 12$.}
      \label{fig:triple-height-momentum}
  \end{figure}}

\section{Conclusion}
Our aim in this paper was to derive a mathematically rigorous
one-dimensional Saint-Venant system, extended to include both
precipitation and infiltration effects. We achieved this by going back
to the original two-dimensional Navier--Stokes equations and adapting
the boundary conditions as appropriate to model these additional
phenomena. The new model \eqref{ModelSL} that we have derived includes
additional momentum source and friction terms in comparison to other
models, which becomes special cases of our system. The friction terms
are obtained naturally from the derivation and their presence is
essential in explaining how the velocity of the water-body interacts
with the additional water coming from either precipitation or runoff;
we demonstrated in \S\ref{subsec:friction-alpha} that, for certain
regimes, this model may yield non-physical solutions. We showed in
Theorem~\ref{THM_ESV} that the existence
of these additional terms leads to a model whose energetic consistency
depends solely on the level of assumed rain-induced friction, denoted
by $\alpha$.

In developing a numerical model, existing approaches such as finite
difference or finite volume can be used, but they fail to ensure
certain properties that we would like our method to have. The
alternative approach we took was to instead use a kinetic formulation,
writing our Saint-Venant system as a single kinetic equation which can
then be solved using a finite volume method to find the original
variables $(h, q)$. To demonstrate the applicability and viability of
our system and associated kinetic scheme, we ran a number of numerical
simulations of our model; in \S\ref{subsec:real-world-experiment}, we
compared the accuracy of our model against a real-world experiment,
while in \S\ref{subsec:cascade} we saw that increasing the value of
$\alpha$ slows down the propagation of the flow; these results were
all in line with our expectations and analysis of the model.

Though our Saint-Venant model goes some way to incorporating
precipitation and infiltration effects, that it only includes one
spatial dimension makes its application to modelling realistic
real-world problems limited in scope, particularly on very large
domains. The techniques and approach that we have used herein can be
readily extended to the two-dimensional system, though some
consideration needs to be given to determine the friction terms and
kinetic scheme for this system; we note, for example, that it is not
clear the kinetic formulation we have derived for the one-dimensional
model can be extended to a second dimension, and thus an entirely new
approach may be required. This work is proposed for future research.
\section{Acknowledgements}
PT and OL were funded by the EPSRC
(Engineering and Physical Sciences Research Council) via an EPSRC-CASE
grant in collaboration with Ambiental Technical
Solutions~\citet{Butler::url:Ambiental}.

OL would like to thank ME for the kind hospitality at the Université de Toulon
in the framework of the ``Professeur Invité'' funding programme.

All authors acknowledge the support of the Marie Sk{\l}odowska-Curie
ITN ``ModCompShock'' and the encouragement and useful remarks of
Charalambos Makridakis, Martin Todd, Donatella Donatelli, Chiara
Simeoni and Alexander Antonarakis.
\bibliographystyle{abbrvnat}

\ifthenelse{\boolean{shownotes}}{\printindex}{}
\end{document}